\documentclass[final,3p,times]{elsarticle}

\usepackage{amssymb}
\usepackage[usenames,dvipsnames]{pstricks}
\usepackage{epsfig}
\usepackage{pst-grad} 
\usepackage{pst-plot} 
\usepackage{algorithmic}

\newtheorem{Theorem}{Theorem}[section]
\newtheorem{Lemma}[Theorem]{Lemma}
\newtheorem{Corollary}[Theorem]{Corollary}
\newtheorem{Proposition}[Theorem]{Proposition}
\newtheorem{Remark}[Theorem]{Remark}
\newtheorem{Example}[Theorem]{Example}

\newtheorem{Definition}[Theorem]{Definition}

\def\demos{\noindent{\bf Proof of Proposition 6.1. }}

\def\QED{\hfill$\Box$}

\def\mP{{\mathcal P}}
\def\N{{\mathbb N}}

\def\Z{{\mathbb Z}}

\def\gcd{{\rm{gcd}}}

\def\mK{{\mathcal K}}

\begin{document}

\begin{frontmatter}

\title{Graphs and complete intersection toric ideals}

\author[ull]{I.~Bermejo\corref{cor}\fnref{fn1}}
\ead{ibermejo@ull.es}
\author[ull]{I.~Garc\'{\i}a-Marco\fnref{fn1}}
\ead{iggarcia@ull.es}
\author[cinv]{E.~Reyes\fnref{fn2}}
\ead{ereyes@math.cinvestav.mx}

\cortext[cor]{Corresponding author. Tel: +34 922318161}

\fntext[fn1]{Partially supported by Ministerio de Ciencia e
Innovaci\'on, Spain (MTM2010-20279-C02-02).}

\fntext[fn2]{Partially supported by SNI, Mexico.}

\address[ull]{Facultad de Matem\'aticas, Universidad de La Laguna, 38200 La Laguna, Tenerife, Canary Islands, Spain}
\address[cinv]{Departamento de Matem\'aticas, Centro de Investigaci\'on y de Estudios
Avanzados del IPN, Apartado Postal 14-740, 0700 M\'exico City, D.F.,
Mexico}

\begin{abstract}
Our purpose is to study the family of simple undirected graphs whose
toric ideal is a complete intersection from both an algorithmic and
a combinatorial point of view. We obtain a polynomial time algorithm
that, given a graph $G$, checks whether its toric ideal $P_G$ is a
complete intersection or not. Whenever $P_G$ is a complete
intersection, the algorithm also returns a minimal set of generators
of $P_G$. Moreover, we prove that if $G$ is a connected graph and
$P_G$ is a complete intersection, then there exist two induced
subgraphs $R$ and $C$ of $G$ such that the vertex set $V(G)$ of $G$
is the disjoint union of $V(R)$ and $V(C)$, where $R$ is a bipartite
ring graph and $C$ is either the empty graph, an odd primitive
cycle, or consists of two odd primitive cycles properly connected.
Finally, if $R$ is $2$-connected and $C$ is connected, we list the
families of graphs whose toric ideals are complete intersection.
\end{abstract}


\end{frontmatter}

\parskip 8pt

\section{Introduction}\label{Int}

Let $k$ be an arbitrary field and $A=(a_{ij})$ an $m \times n$
matrix with non negative integer entries $a_{ij}$ and with non-zero
columns. Let $k[x_1,\ldots,x_n]$ and $k[t_1,\ldots,t_m]$ be two
polynomial rings over $k$. Denote by $x^b$ the monomial $x_1^{b_1}
\cdots x_n^{b_n}$, where $b=(b_1,\ldots,b_n)\in \N^n$. A {\it
binomial} $f$ in $k[x_1,\ldots,x_n]$ is a difference of two
monomials, i.e., $f=x^b-x^c$ for some $b,c\in\N^n$. An ideal
generated by binomials is called a {\it binomial ideal}. Consider
$\varphi$ the graded homomorphism of $k$-algebras
$$\varphi\colon k[x_1,\ldots,x_n] \rightarrow k[t_1,\ldots,t_m] \,
\mbox{ induced by }\, \varphi(x_i)=t^{a_i},
$$
where $a_i$ is the $i${\it-th} column of $A$. The polynomial rings
are graded by assigning $\deg(t_i) = 1$ and
$\deg(x_j)=\deg(t^{a_j})$ for every $i,j$. The kernel of $\varphi$,
denoted by $P_A$, is called the {\it toric ideal\/} associated to
$A$. It is well-known that $P_A$ is a prime graded binomial ideal
with ${\rm ht}(P_A) = n - {\rm rank}(A)$ (see for example
\cite{Stur1, monalg}).

$P_A$ is a {\it complete intersection} if $\mu(P_A) = {\rm
ht}(P_A)$, where $\mu(P_A)$ denotes the minimal number of generators
of $P_A$. Equivalently, $P_A$ is a complete intersection if and only
if there exists a set of homogeneous binomials $f_1,\ldots,f_r \in
k[x_1,\ldots,x_n]$ such that $r = {\rm ht}(P_A)$ and
$P_A=(f_1,\ldots,f_r)$.

Complete intersection toric ideals were first studied by Herzog in
\cite{He3}. After that, they have been extensively studied by
several authors; see for example \cite{stcib-algorithm,
 BGsimplicial, morales-thoma} and the references
there. It is well known, see e.g. \cite{F-M-S-2} or \cite{SSS}, that
the problem of deciding whether a toric ideal is a complete
intersection belongs to the complexity class $\mathcal{NP}$.

Let $G$ be a simple undirected graph, i.e., an undirected graph
without multiple edges or loops. Set $V(G) = \{v_1, \ldots, v_m\}$
its vertex set, $E(G) = \{e_1,\ldots, e_n\}$ its edge set, and $A_G$
its incidence matrix. The toric ideal associated to $A_G$ is denoted
by $P_G$. It is a prime homogeneous binomial ideal called the {\it
toric ideal of $G$}. The image of $\varphi$ is denoted by $k[G]$ and
called the {\it edge algebra of $G$}. If we denote by $b(G)$ the
number of connected components of $G$ which are bipartite, then
${\rm rank}(A_G) = m - b(G)$ (see \cite{Vil2}) which implies that
${\rm ht}(P_G) = n - m + b(G)$. We say that $G$ is a {\it complete
intersection} if the corresponding toric ideal $P_G$ is a complete
intersection.

In this work we study the complete intersection property of graphs
from both an algorithmic and a combinatorial point of view.

The complete intersection property for bipartite graphs has been
extensively studied; see for example \cite{luisa-tor, GRV, Ring,
accota-gv,katzman,aron-jac}. It is worth mentioning that Gitler,
Reyes and Villarreal proved in \cite{Ring} that a bipartite graph is
a complete intersection if and only if it is a ring graph. Since
ring graphs are obviously planar, they could derive that every
complete intersection bipartite graph is planar, which was
previously proved by Katzman \cite{katzman} without using the notion
of ring graph. When graphs are not necessarily bipartite there is
some recent work by Tatakis and Thoma \cite{Tatakis-Thoma}, in the
last section we make use of some of their technical results. For
directed graphs, the complete intersection property has also been
widely studied, see for example \cite{GRVega, Ring, Morfismos}.

\begin{figure}
\begin{center}
\scalebox{1} {
\begin{pspicture}(0,-0.89)(3.42,0.89)
\psdots[dotsize=0.07](0.66,0.07) \psdots[dotsize=0.07](0.6,-0.65)
\psdots[dotsize=0.07](0.08,-0.33)
\psline[linewidth=0.01cm](0.08,-0.33)(0.62,0.07)
\psline[linewidth=0.01cm](0.66,0.11)(0.62,-0.67)
\psline[linewidth=0.01cm](0.62,-0.67)(0.06,-0.33)
\psline[linewidth=0.01cm](0.68,0.11)(1.6,0.09)
\psdots[dotsize=0.07](1.6,0.09) \psdots[dotsize=0.07](1.12,0.79)
\psdots[dotsize=0.07](1.96,0.79) \psdots[dotsize=0.07](2.82,0.07)
\psdots[dotsize=0.07](2.42,-0.79) \psdots[dotsize=0.07](3.32,-0.77)
\psline[linewidth=0.01cm](1.12,0.79)(1.58,0.13)
\psline[linewidth=0.01cm](1.1,0.79)(1.94,0.79)
\psline[linewidth=0.01cm](1.94,0.79)(1.58,0.11)
\psline[linewidth=0.01cm](1.6,0.11)(2.84,0.09)
\psline[linewidth=0.01cm](2.82,0.09)(2.44,-0.77)
\psline[linewidth=0.01cm](2.44,-0.77)(3.34,-0.77)
\psline[linewidth=0.01cm](3.3,-0.77)(2.82,0.11)
\end{pspicture}
 }
 \vskip-.2cm
\caption{Ring graph which is not a complete
intersection}\label{fig1}
\end{center}
\end{figure}
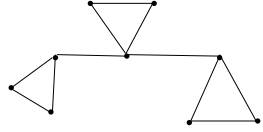

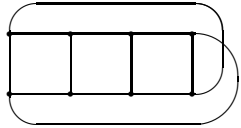
\begin{figure}
\setlength{\unitlength}{.02cm}
\begin{center}
\begin{picture}(0,45)(0,45)

\put(10,70){\oval(140,40)[t]}

\put(60,60){\oval(40,60)[r]}

\put(15,30){\oval(150,40)[lb]}

\put(60,40){\oval(60,60)[r]}

\put(-60,70){\line(0,-1){40}}

\put(-60,70){\line(0,-1){40}}

\put(20,70){\line(0,-1){40}}

\put(-20,70){\line(0,-1){40}}

 \put(-60,70){\line(1,0){80}}

\put(-60,70){\circle*{4}}

\put(-20,70){\circle*{4}}

\put(20,70){\circle*{4}}

\put(-60,30){\line(1,0){80}}

\put(-60,30){\circle*{4}}

\put(-20,30){\circle*{4}}

\put(20,30){\circle*{4}}

\put(20,30){\line(1,0){40}}

\put(20,70){\line(1,0){40}}

\put(60,30){\circle*{4}}

\put(60,70){\circle*{4}}

\put(60,30){\line(0,1){40}}

\put(15,10){\line(1,0){50}}
\end{picture}
\vspace{.5cm} \caption{Non planar graph which is a complete
intersection}\label{fig2}
\end{center}
\end{figure}

In this work, our graphs are undirected and not necessarily
bipartite. In this general setting, the problem requires a different
approach. Indeed, Figure \ref{fig1} shows an example of a ring graph
whose toric ideal is not a complete intersection. Moreover, there
exist complete intersection graphs which are not ring graphs; Figure
\ref{fig2} shows a complete intersection graph which is not even
planar.

The main results of this work are Theorem \ref{principal}, Theorem
\ref{teoremaestructura} and Theorem \ref{ultimo}. The first one
yields a polynomial time algorithm which receives as input a simple
undirected graph $G$ and returns {\sc True} if $G$ is a complete
intersection or {\sc False} other\-wise. Moreover, whenever $G$ is a
complete intersection, the algorithm provides without any extra
effort a minimal set of generators of $P_G$. As a consequence of
this algorithm we obtain that the problem of determining whether a
graph $G$ is a complete intersection belongs to the complexity class
$\mathcal P$. Given a connected graph $G$, we get a partition of $G$
into two disjoint induced subgraphs $C$ and $R$ such that $V(C) =
V(C_1) \bigsqcup \cdots \bigsqcup V(C_s)$ where $C_1,\ldots,C_s$ are
odd primitive cycles, and $R$ is bipartite. In this context, Theorem
\ref{teoremaestructura} gives necessary conditions for a graph to be
a complete intersection by characterizing when $C$ is a complete
intersection. Using this result, when $C$ is connected and $R$ is
$2$-connected, Theorem \ref{ultimo} characterizes the complete
intersection property on $G$ by determining all possible edges
connecting $C$ and $R$.

In Section \ref{sec2}, we collect some results concerning general
toric ideals that will be useful in the sequel. The main result in
this section is Proposition \ref{subconj2}, which deals with the
problem of when the complete intersection property is preserved by
elimination of variables. For toric ideals associated to graphs,
Proposition \ref{subconj2} states that any induced subgraph of a
complete intersection graph also has this property. This is Theorem
\ref{induce-CI} in Section \ref{sec3}, which allows us to obtain in
Theorem \ref{2cases} an upper bound for the number of edges of a
complete intersection graph in terms of the number of vertices,
improving all previously known bounds (see Corollary
\ref{cotasuperior}). An immediate consequence of Theorem
\ref{2cases} is that a complete intersection graph either has a
vertex of degree $\leq 2$, or is $3$-regular (see Corollary
\ref{2situaciones}). Section \ref{sec4} is devoted to designing
Algorithm CI-graph, a polynomial time algorithm for checking whether
a graph is a complete intersection. This algorithm is a direct
consequence of Theorem \ref{principal} and works as follows:
vertices of degree $1$ are removed, also vertices of degree $2$ are
removed after checking certain conditions; if these conditions are
not satisfied, the algorithm returns {\sc False}; otherwise, we
iterate this process until we get either a trivial graph or a graph
in which every vertex has degree $\geq 3$. If there is a vertex of
degree $> 3$, the algorithm returns {\sc False}. Otherwise we use
the characterization of complete intersection $3$-regular graphs
given in Theorem \ref{todosgradomayor2}. Finally, we use Theorem
\ref{fs} to check if $G$ is a complete intersection. Section
\ref{sec5} deals with the problem of finding forbidden subgraphs in
a complete intersection graph. The main result is Theorem
\ref{notheta}, where we prove that odd theta graphs whose base
vertices are not adjacent, and also even theta graphs, are forbidden
subgraphs of a complete intersection graph (see Definition
\ref{thetadefinition} for a definition of even and odd theta
graphs). To prove this, we use Lemma \ref{2subgrafos} and
Proposition \ref{contraccion}, two technical results concerning the
vertices of degree $2$ in a complete intersection graph. In Section
\ref{sec7} we apply the previous results in order to obtain the
above mentioned Theorem \ref{teoremaestructura} and Theorem
\ref{ultimo} together with their normal versions; Corollary
\ref{teoremaestructuranormal} and Corollary \ref{ultimonormal}.

\section{Complete Intersection toric ideals} \label{sec2}

In this section, $A$ denotes an $m \times n$ matrix with non-zero
columns $a_1,\ldots,a_n \in \N^m$ and $P_A \subset
k[x_1,\ldots,x_n]$ is the toric ideal of $A$, which is the kernel of
the $k$-algebra homomorphism $\varphi: k[x_1,\ldots,x_n]
\longrightarrow k[t_1,\ldots,t_m]$ induced by $\varphi(x_i) =
t^{a_i}$.

\begin{Definition}Let $T$ be a subset of $\{t_1,\ldots,t_m\}$. We define
$T_{\varphi^{-1}}$ as the set $\{x_i \, \vert \, \varphi(x_i) \in
k[T]\}$.
\end{Definition}

We have that $P_A \cap k[T_{\varphi^{-1}}]$ is the toric ideal
associated to the matrix whose columns are the $i$-th columns of
 $A$ such that $x_i \in T_{\varphi^{-1}}$; see
\cite[Proposition 4.13(a)]{Stur1}.

\begin{Lemma}\label{subconj1}
Let $\mathfrak B$ be a set of generators of $P_A$ consisting of
binomials, then $\mathfrak B  \cap  k[T_{\varphi^{-1}}]$ is a set of
generators of $P_A  \cap  k[T_{\varphi^{-1}}]$. Moreover, if
$\mathfrak B$ is minimal, then $\mathfrak B \cap
k[T_{\varphi^{-1}}]$ is minimal.
\end{Lemma}
\begin{demo}Our proof begins the observation that whenever $f = x^{\alpha} - x^{\beta} \in
P_A$, then $x^{\alpha} \in k[T_{\varphi^{-1}}]$ if and only if
$x^{\beta} \in k[T_{\varphi^{-1}}]$, by the definition of
$T_{\varphi^{-1}}$. Now let $g$ be a binomial in $P_A \cap
k[T_{\varphi^{-1}}]$. Since $g \in P_A$, then $g = \sum_{f_i \in
\mathfrak B} g_i\, f_i$ with $g_i \in k[x_1,\ldots,x_n]$. Consider
the morphism $\psi$ defined by $\psi(x_i) = x_i$ if $x_i \in
T_{\varphi^{-1}}$, and $\psi(x_i) = 0$ otherwise. Then for every
$f_i \in \mathfrak B$, we get that $\psi(f_i)= f_i$ if $f_i \in
k[T_{\varphi^{-1}}]$ or $\psi(f_i) = 0$ otherwise. Thus, $$g =
\psi(g) = \sum_{f_i \in \mathfrak B} \psi(g_i)\, \psi(f_i) =
\sum_{f_i \in \mathfrak B \cap k[T_{\varphi^{-1}}]} \psi(g_i)\,
f_i.$$ Hence, $\mathfrak B \cap k[T_{\varphi^{-1}}]$ generates $P_A
\cap k[T_{\varphi^{-1}}]$. Moreover, $\mathfrak B \cap
k[T_{\varphi^{-1}}]$ is minimal whenever $\mathfrak B$ is. \QED
\end{demo}

\begin{Proposition}\label{subconj2}If
$P_A$ is a complete intersection, then $P_A \cap
k[T_{\varphi^{-1}}]$ is a complete intersection.
\end{Proposition}

\begin{demo}Let $\mathfrak B$ be a minimal set of generators of $P_A$
consisting of binomials, then $\mathfrak B$ is a regular sequence.
Hence, by  Lemma \ref{subconj1}, the set $\mathfrak B \cap
k[T_{\varphi^{-1}}]$  is a regular sequence which generates $P_A
\cap k[T_{\varphi^{-1}}]$ and the result follows. \QED
\end{demo}

\smallskip
Note that Proposition \ref{subconj2} provides a new proof of
\cite[Theorem 4.1]{KMT} for the particular case of a toric ideal.

In the following sections we will use Theorem 1.1 in \cite{HSh},
which is a reformulation of \cite[Theorem 2.9]{Fisher-Shapiro}. For
presenting this result we have to introduce first some definitions.

Let $B$ be an integral matrix, $B$ is called {\it mixed} if every
row of $B$ has a positive and a negative entry. $B$ is said to be
{\it dominating} if it does not contain any square mixed submatrix.
$\Delta_t(B)$ denotes the greatest common divisor of every $t \times
t$ minor of $B$ where $t \leq {\rm rank}(B)$.

\begin{Theorem}(\cite[Theorem 1.1]{HSh}, \cite[Theorem 2.9]{Fisher-Shapiro})
\label{fs} Let $P_A$ be a toric ideal of height $r$ and $g_i  =
x^{\alpha_i} - x^{\beta_i} \in P_A$ with
 $\gcd (x^{\alpha_i}, x^{\beta_i}) = 1$ for $1 \leq i \leq r$ .
If $B$ denotes the $r \times n$ matrix whose $i$-th row is $\alpha_i
- \beta_i$ for $1 \leq i \leq r$, then
\begin{center} $P_A = (g_1,\ldots,g_r) \,
 \Longleftrightarrow \, B$  is dominating and $\Delta_r(B) =  1$.
 \end{center} \end{Theorem}

The following result, whose proof is straightforward, will be useful
to prove that certain matrices are dominating.

\begin{Lemma}\label{dominating}Let $B$ be an $r \times n$ matrix
with column vectors $c_1,\ldots,c_n \in \Z^r$ such that $c_i$ has
only one nonzero entry for some $1 \leq i \leq n$, and denote by
$B'$ the $r \times n-1$ matrix with column vectors
$c_1,\ldots,c_{i-1},c_{i+1},\ldots,c_n$. Then, $B$ is dominating if
and only if $B'$ is dominating.
\end{Lemma}

\section{An upper bound for the number of edges in a complete intersection graph} \label{sec3}

We begin this section by setting up some notation and terminology
about graphs. For unexplained terminology and results on graphs we
refer to \cite{diestel, Harary}.

A {\it walk} $w$ connecting $u, v \in V(G)$ is a finite sequence of
vertices $w = (u=v_{i_0}, v_{i_1},\ldots, v_{i_q} = v)$ such that
$\{v_{i_{j-1}},v_{i_j}\} \in E(G)$ for every $1 \leq j \leq q$. If
$v_{i_j} \neq v_{i_k}$ for every $0 \leq j < k \leq q$ then $w$ is
called a {\it path}. The vertex set of the walk $w$ is $V(w) :=
\{v_{i_0},\ldots,v_{i_{q}}\}$ and its edge set is $E(w) :=
\{\{v_{i_{j-1}},v_{i_j}\}\,\vert\, 1 \leq j \leq q\}$. The {\it
length of the walk} is the number $q$ of edges in the walk. An {\it
even} (respectively {\it odd}) {\it walk} is a walk of even
(respectively odd) length. A walk is {\it closed} if $u = v$. A {\it
cycle} is a closed walk with $v_{i_k} \neq v_{i_j}$ for every $1
\leq k < j \leq q$. A cycle is {\it primitive} if
$\{v_{i_k},v_{i_j}\} \notin E(G)$ for every $1 \leq k < k+1 < j \leq
q$.

For a walk $w = (u=v_{i_0}, v_{i_1},\ldots, v_{i_q} = v)$ we denote
by $-w$ the {\it inverse walk} $(v =
v_{i_{q}},\ldots,\,v_{i_1},\,v_{i_0} = u)$. Let $w_1,\ldots,w_r$ be
walks such that $w_i$ connects $u_i, u_{i+1}$ for every $i \in
\{1,\ldots,r\}$, then $(w_1,\ldots,w_r)$ denotes the walk connecting
$u_1, u_{r+1}$ obtained by sticking the walks $w_1,\ldots,\,w_{r-1}$
and $w_r$ together.

 Given an even closed walk, $w =
(v_{i_0},\ldots, v_{i_{2q}} = v_{i_0})$ where $e_{k_j} =
\{v_{i_{j-1}}, v_{i_j}\}$ for $1 \leq j \leq 2q$, we denote by $B_w$
the binomial
$$B_w := \prod_{l = 1}^q x_{k_{2l-1}}- \prod_{l = 1}^q x_{k_{2l}}.$$

Villarreal  \cite[Proposition 3.1]{Vi3} proved that $P_G$ is
generated by these binomials, i.e., $P_G = (\{B_w \, \vert \, w$ is
an even closed walk$\})$. Hibi and Ohsugi \cite[Lemma 3.2]{OH2}
improved this result by giving a necessary condition for a binomial
in $P_G$ to be primitive. Recall that $x^{\alpha} - x^{\beta} \in
P_G$ is {\it primitive} if there exists no other binomial
$x^{\alpha'} - x^{\beta'} \in P_G$ such that $x^{\alpha'} \mid
x^{\alpha}$ and $x^{\beta'} \mid x^{\beta}$. Whenever a binomial
belongs to a minimal set of generators of $P_G$, then it is
necessarily primitive (see \cite{Stur1}); thus the set of all
primitive binomials of $P_G$, which is called the {\it Graver basis
of} $P_G$, is a set of generators for $P_G$.

\begin{Lemma} \cite[Lemma 3.2]{OH2} \label{primitivo} If $B_w$ is primitive, then
one of these holds:
\begin{itemize}
\item $w$ is an even cycle,
\item $w = (C_1, C_2)$ where $C_1$ and $C_2$ are odd cycles
having exactly a vertex in common, or
\item $w = (C_1, w_1, C_2, -w_2)$ where $C_1, C_2$ are
vertex disjoint odd cycles and $w_1, w_2$ are walks connecting a
vertex $v_1 \in V(C_1)$ and a vertex $v_2 \in V(C_2)$.
\end{itemize}
\end{Lemma}

 For a complete characterization of primitive binomials and a description of all
minimal sets of generators of $P_G$ formed by binomials we refer the
reader to \cite{Reyes-Tatakis-Thoma}.

Now we aim to prove that the complete intersection is hereditary,
i.e., if a graph is a complete intersection then every induced
subgraph also is. Let us first recall the definition of induced
subgraph.

\begin{Definition}Let $G$ be a graph, $G'$ is an {\it induced subgraph of $G$} if
$V(G') \subset V(G)$ and \begin{center} $E(G') = \{e \in E(G)\,
\vert\, e \subset V(G')\}$.
\end{center} If $V' \subset V(G)$, we will denote by $[V']$ the induced subgraph of $G$ with
vertex set $V'$. Let $v_1,\ldots,v_s$ be vertices of $G$, the
induced subgraph  $[V(G) \setminus \{v_1,\ldots,v_s\}]$ will also be
denoted by $G \setminus \{v_1,\ldots,v_s\}.$
\end{Definition}

For an induced subgraph $G'$, if we denote $T := \{t_i \, \vert \,
v_i \in V(G')\}$, then $T_{\varphi^{-1}} = \{x_i  \, \vert\, e_i \in
E(G')\}$ and $P_{G'}=P_G \cap k[T_{\varphi^{-1}}]$. Hence, by Lemma
\ref{subconj1} and Proposition \ref{subconj2}  we deduce the
following results.

\begin{Proposition} \label{generadoresinducido}
Let $G'$ be an induced subgraph of $G$. If $P_G =
(B_{w_1},\ldots,B_{w_s})$ for some even closed walks
$w_1,\ldots,w_s$ in $G$, then $P_{G'} = (B_{w_i}\,\vert\, V(w_i)
\subset V(G'),\, 1 \leq i \leq s)$.
\end{Proposition}

\begin{Theorem}\label{induce-CI}
Let $G'$ be an induced subgraph of $G$. If $G$ is a complete
intersection, then so is $G'$.
\end{Theorem}

A different proof of Theorem \ref{induce-CI} exists also in
\cite[Theorem 3.1]{Tatakis-Thoma}.

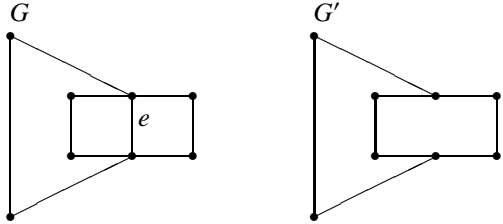
\begin{figure}[h]
\begin{center}
\setlength{\unitlength}{.04cm}
\begin{picture}(150,65)(20,10)

\put(20,75){$G$}

\put(20,70){\circle*{3}}

\put(20,70){\line(0,-1){60}}

\put(20,70){\line(2,-1){40}}

\put(20,10){\circle*{3}}

\put(20,10){\line(2,1){40}}

\put(40,50){\circle*{3}}

\put(40,50){\line(1,0){40}}

\put(62,40){$e$}

\put(60,50){\line(0,-1){20}}

\put(60,50){\circle*{3}}

\put(80,50){\circle*{3}}

\put(40,30){\circle*{3}}

\put(40,30){\line(1,0){40}}

\put(40,30){\line(0,1){20}}

\put(80,30){\line(0,1){20}}

\put(60,30){\circle*{3}}

\put(80,30){\circle*{3}}


\put(120,75){$G'$}

 \put(120,70){\circle*{3}}

\put(120,70){\line(0,-1){60}}

\put(120,70){\line(2,-1){40}}

\put(120,10){\circle*{3}}

\put(120,10){\line(2,1){40}}

\put(140,50){\circle*{3}}

\put(140,50){\line(1,0){40}}

\put(160,50){\circle*{3}}

\put(180,50){\circle*{3}}

\put(140,30){\circle*{3}}

\put(140,30){\line(1,0){40}}

\put(140,30){\line(0,1){20}}

\put(180,30){\line(0,1){20}}

\put(160,30){\circle*{3}}

\put(180,30){\circle*{3}}
\end{picture}
\caption{$G'$ is a subgraph of $G$, both are bipartite but $G$ is a
ring graph and $G'$ is not. Thus $G$ is a complete intersection and
$G'$ is not.}\label{fig3}
\end{center}
\end{figure}

\smallskip These results are not true in general if we drop the
assumption that $G'$ is induced, as the example in Figure \ref{fig3}
shows.

An almost immediate consequence of Theorem \ref{induce-CI} is that a
graph is a complete intersection if and only if all its connected
components are complete intersections. This allows us to reduce our
study to connected graphs.

\begin{Corollary}\label{conexo}Let $G$ be a graph with connected components
$G_1,\ldots,G_s$. Then, $G$ is a complete intersection if and only
if so are $G_1,\ldots,G_s$.
\end{Corollary}
\begin{demo}One implication follows by Theorem \ref{induce-CI},
because $G_i$ is an induced subgraph of $G$ for $1 \leq i \leq s$.
Assume now that $G_1,\ldots,G_s$ are complete intersections and let
$\mathfrak B_i$ be a minimal set of generators of $P_{G_i}$ for $1
\leq i \leq s$. Since $P_G = (\{B_w \,\vert \, w$ is an even closed
walk$\})$ and every even closed walk is necessarily contained in a
connected component of $G$, it is evident that $P_G = \langle
\mathfrak B_1 \cup \cdots \cup \mathfrak B_s \rangle$. Moreover, we
have that ${\rm ht}(P_G) = {\rm ht}(P_{G_1}) + \cdots + {\rm
ht}(P_{G_s})$, and the result follows. \QED
\end{demo}

For a vertex $v \in V(G)$, the {\it neighborhood of $v$} is the set
of vertices which are adjacent to $v$, i.e.,  $N_G(v) := \{u \in
V(G)\, \vert \, \{u,v\} \in E(G) \}$. The cardinality of this set is
called the {\it degree of} $v$ and is denoted by ${\rm deg}_G(v)$,
or ${\rm deg}(v)$ when $G$ is understood. $G$ is $k$-regular if
every vertex of $G$ has degree $k$.

The following result provides an upper bound for the number of edges
of complete intersection graphs. This gives the taste that they can
not be very dense.

\begin{Theorem}\label{2cases}
Let $G$ be a complete intersection connected graph, then:
\begin{itemize}
\item $2\, |E(G)| + 4 \leq 4\, |V(G)| - \sum_{v \in V(G)} b(G \setminus \{v\})$ if $G$ is bipartite.
\item $2\, |E(G)| \leq 3\, |V(G)| - \sum_{v \in V(G)} b(G \setminus \{v\})$ if $G$ is not bipartite.
\end{itemize}
In both cases equality is attained if and only if $P_G$ is generated
by quadrics.
\end{Theorem}
\begin{demo}Let $\{B_{w_1},\ldots,B_{w_r}\}$ be a minimal set of generators of
$P_G$, where $w_i$ is an even closed walk for every $1 \leq i \leq
r$. By Proposition \ref{generadoresinducido}, for every $v \in V(G)$
we have that $P_{G \setminus \{v\}} = ( B_{w_i} \,\vert\, v \not\in
V(w_i),\, 1 \leq i \leq r)$ and by Lemma \ref{primitivo}, it follows
that $\vert V(w_i) \vert \geq 4$ and $\vert V(w_i) \vert = 4$ if and
only if $w_i$ is a cycle of length $4$, which is equivalent to
$B_{w_i}$ is a quadric.

Therefore,
$$4 \mu(P_G) \leq \sum_{v\in V(G)} \mu(P_G) - \mu(P_{G \setminus
\{v\}})$$ and equality holds if and only if $P_G$ is generated by
quadrics. Now suppose that $G$ is a complete intersection then, by
Theorem \ref{induce-CI}, $G \setminus \{v\}$ is a complete
intersection for every $v \in V(G)$. Hence, $$4 {\rm ht}(P_G) \leq
\sum_{v\in V(G)} \left( \mu(P_G) - \mu(P_{G \setminus \{v\}})
\right) = \sum_{v\in V(G)} \left({\rm ht}(P_G) - {\rm
ht}(P_{G\setminus \{v\}} )\right) = $$
$$= \sum_{v\in V(G)} ({\rm deg}(v) - 1 + b(G) - b(G \setminus
\{v\})) = 2n - m + b(G)\, m - \sum_{v\in V(G)} b(G \setminus
\{v\}).$$ If $G$ is bipartite, then $b(G) = 1,\, {\rm ht}(P_G) = n -
m + 1$ and $2n + 4 \leq 4m - \sum_{v\in V(G)} b(G \setminus \{v\}),$
and if $G$ is not bipartite, then $b(G) = 0,\, {\rm ht}(P_G) = n -
m$ and  $2n \leq 3m - \sum_{v\in V(G)} b(G \setminus \{v\}).$
 In both cases equality is attained if and only if
$P_G$ is generated by quadrics.\QED
\end{demo}

\bigskip

Katzman in {\rm \cite[Corollary 3.8]{katzman}} proved that $|E(G)| +
4 \leq 2 \, |V(G)|$ for a complete intersection connected bipartite
graph. Independently, from a result of {\rm Fischer, Morris} and
{\rm Shapiro} {\rm \cite[Corollary 3.4]{F-M-S-2}} one can deduce
that if $G$ is a complete intersection connected graph, then $|E(G)|
+ 4 \leq 2 \, |V(G)|$ if $G$ is bipartite and $|E(G)| + 2 \leq 2 \,
|V(G)|$ if $G$ is non bipartite. The following result improves these
bounds.

\begin{Corollary}\label{cotasuperior}Let $G$ be a complete intersection connected graph, then
\begin{itemize}
\item $2\, |E(G)| + 4 \leq 3\, |V(G)|$ if $G$ is bipartite, and
\item $2\, |E(G)| \leq 3\, |V(G)|$ if $G$ is not bipartite.
\end{itemize}
\end{Corollary}
\begin{demo}It is a consequence of Theorem \ref{2cases} and that if
$G$ is bipartite, then $b(G \setminus \{v\}) \geq 1$ for every $v
\in V(G)$. \QED
\end{demo}

This section ends with two more consequences of Theorem
\ref{2cases}.

\begin{Corollary}\label{2situaciones}
Let $G$ be a complete intersection graph, then either
\begin{itemize}[leftmargin=1cm]
\item[{\rm (a)}] there exists a vertex of degree $\leq 2$, or
\item[{\rm (b)}] $G$ is $3$-regular, $b(G \setminus \{v\}) = 0$ for every $v \in V(G)$
and $P_G$ is generated by quadrics.
\end{itemize}
\end{Corollary}
\begin{demo}Assume that $G$ is a connected graph and every vertex has degree $\geq
3$. Then $2n = \sum_{v \in V(G)}{\rm deg}(v) \geq 3m$; hence by
Corollary \ref{cotasuperior}, $G$ is not bipartite and $2n = 3m$.
Thus, $G$ is $3$-regular and by Theorem \ref{2cases} this can only
happen if (b) holds. \QED
\end{demo}

\bigskip We denote by $\mK_m$ the complete graph with $m$ vertices
and by $\mK_{m_1,m_2}$ the complete bipartite graph with partitions
of sizes $m_1$ and $m_2$.

\begin{Corollary}\label{subgraphK23}
If $G$ is a complete intersection, then it does not contain
$\mK_{2,3}$ as a subgraph.
\end{Corollary}
\begin{demo}Assume that $G$ contains $\mK_{2,3}$ as a subgraph and denote by $H$
the induced subgraph of $G$ with $5$ vertices containing $\mK_{2,3}$
as a subgraph. If $H = \mK_{2,3}$ then $2 |E(H)| + 4 = 16
> 15 = 3 |V(H)|$. If $E(H) = E(\mK_{2,3}) \cup \{e_1,\ldots,e_s\}$,  we
have that $H$ is not bipartite and if $s = 1$ and $e_1 =
\{v_1,\,v_2\} \subset V(H)$, then $b(H \setminus \{v_i\}) \geq 1$
for $i = 1,2$. So
$$2 \, \vert E(H) \vert \, = 12 + 2s  >  15 + 2(s-2) \geq 3 \, \vert V(H) \vert  - \sum_{v \in V(H)}
b(H \setminus \{v\}).$$ In both cases, one gets that $H$ is not a
complete intersection by Theorem \ref{2cases} and Corollary
\ref{cotasuperior}. Furthermore, by Theorem \ref{induce-CI} one
concludes that $G$ is not a complete intersection. \QED
\end{demo}

\section{The algorithm} \label{sec4}

The aim of this section is to provide {\rm Algorithm CI-graph}, an
algorithm for checking whether a graph is a complete intersection.
This algorithm follows as a consequence of Theorem \ref{principal},
which is the main result of this section.

By Corollary \ref{2situaciones} we have that a complete intersection
graph either has a vertex of degree $\leq 2$ or is $3$-regular. This
section begins with a thorough study of $3$-regular complete
intersection graphs. It will turn out in Theorem
\ref{todosgradomayor2} that a $3$-regular graph is a complete
intersection if and only if it is an odd band or an even M\"obius
band. To prove this we need some definitions and a technical lemma.
Theorem \ref{todosgradomayor2} will be essential for proving Theorem
\ref{principal}.

\begin{Definition}
A {\it chain} is a graph $G$ with $V(G)=\{a_1,\ldots, a_r,
b_1,\ldots, b_r\}$ and edges $\{a_i, a_{i+1}\}$,$\{b_i, b_{i+1}\}$
and $\{a_j, b_j\}$ for $1 \leq i < r,\, 1 \leq j \leq r$.
\end{Definition}

\begin{Definition}
Let $G$ be a graph with a subgraph $H$ such that $V(G) = V(H)$ and
$H$ is a chain. If $E(G) = E(H)\, \cup \, \{\{a_1, a_r\},\{b_1,
b_r\}\}$ we say that $G$ is a {\it band}. If $E(G) = E(H)\, \cup  \,
\{\{a_1, b_r\},\{a_r, b_1\}\}$ we say that $G$ is a {\it
M\"obius-band}. In addition if $r$ is odd, we say that $G$ is an
{\it odd M\"obius-band} $($or {\it odd band}$)$ and if $r$ is even,
we say that $G$ is an {\it even M\"obius-band} $(${\it or even
band}$)$.
\end{Definition}

\begin{Lemma}\label{3vecinos}Let $G$ be a complete intersection connected
$3$-regular graph. Then, either $G = \mK_4$ or for every $v \in
V(G)$ there exists a chain subgraph $H$ of $G$ with $6$ vertices,
such that $v \in V(H)$ and ${\rm deg}_H(v) = 3$.
\end{Lemma}
\begin{demo}Firstly note that $P_G$ is generated by quadrics and $b(G \setminus \{v\}) = 0$ for
every $v \in V(G)$ by Corollary \ref{2situaciones}; in particular,
$G$ is not bipartite. Let $\mathfrak B :=
\{B_{w_1},\ldots,B_{w_r}\}$ be a minimal set of generators of $P_G$
where $w_i$ is a length four cycle for $1 \leq i \leq r$. Take $v
\in V(G)$ and denote by $u_1, u_2, u_3$ its neighbors. From one
hand, we have that $${\rm ht}(P_G) - {\rm ht}(P_{G \setminus \{v\}})
= {\rm deg}(v) - 1 + b(G) - b(G \setminus \{v\}) = 2,$$ and by
Proposition \ref{generadoresinducido} and Theorem \ref{induce-CI},
$P_{G \setminus \{v\}}$ is a complete intersection minimally
generated by $\{B_{w_i} \, \vert \, v \notin V(w_i)\}$; thus $|
\{w_i \, \vert\, v\in V(w_i),\, 1 \leq i \leq r\} | = 2$, and we can
assume that $v \in V(w_1) \cap V(w_2)$. From the other hand,
$|N_G(u_i) \cap N_G(u_j)| \leq 2$ for every $1 \leq i < j \leq 3$;
otherwise $\mK_{2,3}$ is a subgraph of $G$, which is impossible by
Corollary \ref{subgraphK23}. Thus we can assume that $w_1 = (v, u_1,
v_1, u_2, v)$ and $w_2 = (v, u_2, v_2, u_3, v)$ for some $v_1, v_2
\in V(G).$ Since $w_1$ and $w_2$ are length $4$ cycles and
$\mK_{2,3}$ is not a subgraph of $G$, we see that $v_1 \neq v_2$.

If $v_1 = u_3$ or $v_2 = u_1$, then $G = \mK_4$. Otherwise there is
a chain subgraph $H$ of $G$ with $V(H) = \{v, u_1, u_2, u_3, v_1,
v_2\}$ and ${\rm deg}_H(v) = 3$. \QED
\end{demo}

\bigskip

\begin{Theorem}\label{todosgradomayor2}Let $G$ be a $3$-regular graph.
Then, $G$ is a complete intersection if and only if the connected
components of $G$ are odd bands or even M\"obius bands.
\end{Theorem}
\begin{demo}By Corollary \ref{conexo}, we can assume that $G$ is
connected.

$(\Rightarrow)$ Since $\mK_4$ is an even M\"obius band we will
assume $G \neq \mK_4$. By Corollary \ref{2situaciones} we get that
$G$ is a $3$-regular graph with $b(G \setminus \{v\}) = 0$ for every
$v \in V(G)$ and by Lemma \ref{3vecinos} there exists a subgraph
$H'$ of $G$ which is a chain with $6$ vertices. Choose $H$ a chain
subgraph of $G$ maximal with respect to $\vert V(H)\vert$; then
$V(H)=\{a_1,\ldots,a_r,b_1,\ldots,b_r\}$ for some $r \geq 3$.
Applying Lemma \ref{3vecinos} with $v = a_r$, we get that there
exist $a_{r+1}, b_{r+1} \in V(G)$ such that $\{a_r,
a_{r+1}\},\,\{b_r, b_{r+1}\}$ and $\{a_{r+1}, b_{r+1}\} \in E(G)$.
We will prove that $\{a_1, b_1\} = \{a_{r+1}, b_{r+1}\}$.

By the maximality of $H$ either $a_{r+1}$ or $b_{r+1}$ belong to
$V(H)$. We assume that $a_{r+1} \in V(H)$, then
 $a_{r+1} = a_1$ or $a_{r+1} = b_1$. If $a_{r+1} = a_1$, by Lemma \ref{3vecinos},
 one can conclude that $b_{r+1} = b_1$
 because  $\{a_1, b_{r+1}\},\, \{b_r, b_{r+1}\}
\in E(G)$ and $a_2 \neq a_r$. If $a_{r+1} = b_1$, since $\{b_1,
b_{r+1}\},\, \{b_r, b_{r+1}\} \in E(G)$, one can conclude that
either $b_{r+1} = a_1$ or $b_{r+1} = b_2$ and $r = 3$. Furthermore,
if $b_{r+1} = b_2$ and $r = 3$ then there is a $\mK_{2,3}$ subgraph
with vertices $a_1, a_2, a_3, b_1, b_2$ and this is not possible by
Corollary \ref{subgraphK23}.

Therefore we have proved that $H$ is either a band or a M\"obius
band and $G$ is $3$-regular and connected, then $G = H$. Finally,
$G$ can be neither an even band nor an odd M\"obius band, because
both are bipartite and by Corollary \ref{cotasuperior} $G$ is not
bipartite.

$(\Leftarrow)$ Denote $e_i := \{a_i, b_i\}$ for $1 \leq i \leq r$,\,
$e_{r+i} := \{a_i, a_{i+1}\}$ and $e_{2r + i} := \{b_i, b_{i+1}\}$
for $1 \leq i < r$. If $G$ is an odd band, we set $e_{2r} := \{a_1,
a_r\}$ and $e_{3r} := \{b_1, b_r\}$ and if $G$ is an even M\"obius
band, we set $e_{2r} := \{a_1, b_r\}$ and $e_{3r} := \{a_r, b_1\}$.
In both cases $G$ is not bipartite, furthermore $G$ has $3r$ edges
and $2r$ vertices, then ${\rm ht}(P_G) = r$.

Let $w_i$ be the length 4 cycle $w_i := (a_i, b_i, b_{i+1}, a_{i+1},
a_i)$ for $1 \leq i \leq r-1$. If $G$ is an odd band, we denote $w_r
:= (a_1, b_1, b_r, a_r, a_1)$ and if $G$ is an even M\"obius band we
denote $w_r := (a_1, b_1, a_r, b_r, a_1)$. In both cases we have
that $B_{w_i} = x_i  x_{i+1} - x_{r+i} x_{2r+i}$ for $1 \leq i < r$
and $B_{w_r} = x_1  x_r - x_{2r}  x_{3r}$. We denote by
$\{e_1,\ldots,e_{3r}\}$ the canonical basis of $\Z^r$, $\gamma_i :=
e_i + e_{i+1} - e_{r+i} - e_{2r+i}$ for $1 \leq i < r$, $\gamma_r =
e_1 + e_r - e_{2r} - e_{3r}$ and $B$ the $r \times 3r$ matrix whose
rows are $\gamma_1,\ldots,\gamma_r$. Then $\Delta_r(B) = 1$. Let
$B'$ be the $r \times r$ submatrix of $B$ consisting of its
 first $r$ columns. Since every entry of $B'$ is nonnegative we get that
$B'$ is dominating. Furthermore, for every $j > r$ the $j$-th column
of $B$ has only one nonzero entry whose value is $-1$, then by Lemma
\ref{dominating} B is dominating. By Theorem \ref{fs} we can
conclude that $G$ is a complete intersection. \QED
\end{demo}

\bigskip The proof above gives more, whenever $G$ is an odd band or
an even M\"obius band we have obtained a minimal set of generators
of the ideal.

\begin{Corollary}Let $G$ be an odd band,
 then $P_G = (B_{w_1},\ldots,B_{w_r}),$ where $w_i := (a_i, b_i, b_{i+1},
a_{i+1}, a_i)$ for $1 \leq i \leq r-1$ and $w_r := (a_1, b_1, b_r,
a_r, a_1)$.
\end{Corollary}

\begin{Corollary}Let $G$ be an even M\"obius band,
 then
$P_G = (B_{w_1},\ldots,B_{w_r})$, where $w_i := (a_i, b_i, b_{i+1},
a_{i+1}, a_i)$ for $1 \leq i \leq r-1$ and $w_r := (a_1, b_1, a_r,
b_r, a_1)$.
\end{Corollary}

\begin{figure}
\begin{center}
\setlength{\unitlength}{.03cm}
\begin{picture}(-280,100)(185,0)

\put(-50,70){\oval(160,40)[t]}

\put(-50,30){\oval(160,40)[b]}

\put(-130,70){\circle*{4}}

\put(-90,70){\line(0,-1){40}}

\put(-10,70){\line(0,-1){40}}

\put(-50,70){\line(0,-1){40}}

\put(-130,70){\line(1,0){120}}

\put(-130,70){\line(0,-1){40}}

\put(-90,70){\circle*{4}}

\put(-50,70){\circle*{4}}

\put(-10,70){\circle*{4}}

\put(-130,30){\circle*{4}}

\put(-130,30){\line(1,0){120}}

\put(-90,30){\circle*{4}}

\put(-50,30){\circle*{4}}

\put(-10,30){\circle*{4}}

\put(-10,30){\line(1,0){40}}

\put(-10,70){\line(1,0){40}}

\put(30,30){\circle*{4}}

\put(30,70){\circle*{4}}

\put(30,30){\line(0,1){40}}

\put(-47,98){$w_5$}

\put(-115,45){$w_1$}

\put(-75,45){$w_2$}

\put(-35,45){$w_3$}

\put(5,45){$w_4$}

\thicklines

\put(-50,70){\oval(170,46)[t]}

\put(-50,33){\oval(170,50)[b]}

\put(-45,93){\vector(-1,0){10}}

\put(-55,8){\vector(1,0){10}}

\put(-135,70){\vector(0,-1){20}}

\put(-135,30){\line(0,1){40}}

\put(35,30){\line(0,1){40}}

\put(35,30){\vector(0,1){20}}

\put(-127,33){\line(1,0){34}}

\put(-127,33){\vector(1,0){17}}

\put(-93,33){\line(0,1){34}}

\put(-93,33){\vector(0,1){17}}

\put(-93,67){\line(-1,0){34}}

\put(-93,67){\vector(-1,0){17}}

\put(-127,67){\line(0,-1){34}}

\put(-127,67){\vector(0,-1){17}}

\put(-87,33){\line(1,0){34}}

\put(-87,33){\vector(1,0){17}}

\put(-53,33){\line(0,1){34}}

\put(-53,33){\vector(0,1){17}}

\put(-53,67){\line(-1,0){34}}

\put(-53,67){\vector(-1,0){17}}

\put(-87,67){\line(0,-1){34}}

\put(-87,67){\vector(0,-1){17}}

\put(-47,33){\line(1,0){34}}

\put(-47,33){\vector(1,0){17}}

\put(-13,33){\line(0,1){34}}

\put(-13,33){\vector(0,1){17}}

\put(-13,67){\line(-1,0){34}}

\put(-13,67){\vector(-1,0){17}}

\put(-47,67){\line(0,-1){34}}

\put(-47,67){\vector(0,-1){17}}

\put(-7,33){\line(1,0){34}}

\put(-7,33){\vector(1,0){17}}

\put(27,33){\line(0,1){34}}

\put(27,33){\vector(0,1){17}}

\put(27,67){\line(-1,0){34}}

\put(27,67){\vector(-1,0){17}}

\put(-7,67){\line(0,-1){34}}

\put(-7,67){\vector(0,-1){17}}


\put(130,70){\oval(140,40)[t]}

\put(180,60){\oval(40,60)[r]}

\put(135,30){\oval(150,40)[lb]}

\put(180,40){\oval(60,60)[r]}

\put(60,70){\line(0,-1){40}}

\put(60,70){\line(0,-1){40}}

\put(140,70){\line(0,-1){40}}

\put(100,70){\line(0,-1){40}}

 \put(60,70){\line(1,0){80}}

\put(60,70){\circle*{4}}

\put(100,70){\circle*{4}}

\put(140,70){\circle*{4}}

\put(60,30){\line(1,0){80}}

\put(60,30){\circle*{4}}

\put(100,30){\circle*{4}}

\put(140,30){\circle*{4}}

\put(140,30){\line(1,0){40}}

\put(140,70){\line(1,0){40}}

\put(180,30){\circle*{4}}

\put(180,70){\circle*{4}}

\put(180,30){\line(0,1){40}}

\put(135,10){\line(1,0){50}}

{\thicklines

\put(75,45){$w_1$}

\put(115,45){$w_2$}

\put(155,45){$w_3$}

 \put(135,8){\line(1,0){50}}

\put(115,8){\vector(1,0){25}}

\put(165,92){\vector(-1,0){25}}

\put(150,96){$w_4$}

\put(129,70){\oval(144,44)[t]}

\put(182,60){\oval(38,64)[r]}

\put(134,30){\oval(154,44)[lb]}

\put(182,40){\oval(58,64)[r]}

\put(57,70){\line(0,-1){40}}

\put(57,70){\vector(0,-1){20}}

\put(183,28){\line(0,1){44}}

\put(183,72){\vector(0,-1){22}}

\put(63,33){\line(1,0){34}}

\put(63,33){\vector(1,0){17}}

\put(97,33){\line(0,1){34}}

\put(97,33){\vector(0,1){17}}

\put(97,67){\line(-1,0){34}}

\put(97,67){\vector(-1,0){17}}

\put(63,67){\line(0,-1){34}}

\put(63,67){\vector(0,-1){17}}

\put(103,33){\line(1,0){34}}

\put(103,33){\vector(1,0){17}}

\put(137,33){\line(0,1){34}}

\put(137,33){\vector(0,1){17}}

\put(137,67){\line(-1,0){34}}

\put(137,67){\vector(-1,0){17}}

\put(103,67){\line(0,-1){34}}

\put(103,67){\vector(0,-1){17}}

\put(143,33){\line(1,0){34}}

\put(143,33){\vector(1,0){17}}

\put(177,33){\line(0,1){34}}

\put(177,33){\vector(0,1){17}}

\put(177,67){\line(-1,0){34}}

\put(177,67){\vector(-1,0){17}}

\put(143,67){\line(0,-1){34}}

\put(143,67){\vector(0,-1){17}}

\thinlines }

\end{picture}
\end{center}\caption{An odd band, an even M\"obius band and the even closed walks corresponding
to a minimal set of generators of each.} \label{evenMobiusband}
\end{figure}
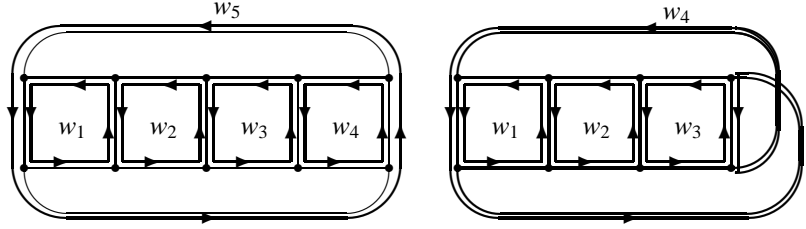

\begin{Remark}Since every even M\"obius band except $\mK_4$ is not
planar, {\rm Theorem \ref{todosgradomayor2}} provides an infinite
family of non planar complete intersection graphs. Both {\rm
Katzman} {\rm \cite{katzman}} and {\rm Gitler, Reyes and Villarreal}
{\rm \cite{Ring}} proved that whenever $G$ is a bipartite complete
intersection then it is planar. As one can see this result is no
longer true if we drop the assumption that $G$ is bipartite. This
was first realized by {\rm Katzman} {\rm \cite[Remark
3.9]{katzman}}, who provided a M\"obius band with $8$ vertices as an
example of a complete intersection non planar graph. Later {\rm
Tatakis and Thoma} {\rm \cite{Tatakis-Thoma}} provided another
example which is not a M\"obius band.
\end{Remark}

We are thus led to the main result of this section.

\begin{Theorem}\label{principal}Let $G$ be a graph without isolated vertices. Then,
$G$ is a complete intersection if and only if one of the following
holds
\begin{enumerate}
\item $\exists\, v \in V(G)$ of degree $1$ and $G \setminus
\{v\}$ is a complete intersection.
\item $\exists\, v \in V(G)$ of degree $2$ such that $b(G \setminus \{v\}) =
b(G) + 1$ and ${G \setminus \{v\}}$ is a complete intersection.
\item $\exists\, v \in V(G)$ of degree $2$ such that $b(G \setminus \{v\}) =
b(G)$, $G \setminus \{v\}$ is a complete intersection and exists a
shortest even closed walk $w$ with
\begin{center}$V(w) = \{v\} \cup N_G(v) \cup \{u \in V(G) \, \vert
\, b(G \setminus \{u,v\}) > b(G \setminus \{u\})\},$ \end{center}
such that
\begin{center} $P_G = P_{G \setminus \{v\}} \cdot k[x_1,\ldots,x_n] +
(B_w).$\end{center}

\item The connected
components of $G$ are odd bands or even M\"obius bands.
\end{enumerate}
\end{Theorem}
\begin{demo}Our proof starts with the observation that if $v \in V(G)$
has degree $1$, then $b(G \setminus \{v\}) = b(G)$ and if it has
degree $2$, then $b(G \setminus \{v\}) - b(G) \in \{0,1\}$. Thus,
the proof falls naturally in the following four cases:

$\begin{array}{rl}$(a)$&$ there exists $v \in V(G)$ such that ${\rm
deg}(v) = 1, \\$(b)$&$ there exists $v \in V(G)$ such that ${\rm
deg}(v) = 2$ and $b(G \setminus \{v\}) = b(G) + 1, \\$(c)$&$ there
exists $v \in V(G)$ such that ${\rm deg}(v) = 2$ and $b(G \setminus
\{v\}) = b(G)$ or $\\$(d)$&$ ${\rm deg}(v) > 2$ for every $v \in
V(G).
\end{array}$

We observe that $J := P_{G \setminus \{v\}} \cdot k[x_1,\ldots,x_n]$
is a prime ideal and $J \subset P_G$. If (a) or (b) holds, then
${\rm ht}(P_{G \setminus \{v\}}) = {\rm ht}(P_G)$, this yields $P_G
= J$ and $G$ is a complete intersection if and only if $G \setminus
\{v\}$ is a complete intersection.

If (c) holds, then ${\rm ht}(P_G) = {\rm ht}(P_{G \setminus \{v\}})
+ 1$. If $G \setminus \{v\}$ is a complete intersection and there
exists an even closed walk $w$ in $G$ such that $P_G = J + (B_w)$,
then $G$ is evidently a complete intersection. Suppose that $G$ is a
complete intersection and let $w_1,\ldots,w_r$ be even closed walks
in $G$ such that $P_G = (B_{w_1},\ldots,B_{w_r})$ with $r = {\rm
ht}(P_G)$. By Proposition \ref{generadoresinducido} and Theorem
\ref{induce-CI}, we have that $P_{G \setminus \{v\}}$ is a complete
intersection minimally generated by $\{B_{w_i} \, \vert \, v \notin
V(w_i)\}$. Since ${\rm ht}(P_{G \setminus \{v\}}) = r - 1$, there
exists a unique $i \in \{1,\ldots,r\}$ such that $v \in V(w_i)$ and
$P_G = J + (B_{w_i}).$ It is obvious that $N_G(v) := \{v_1,v_2\}
\subset V(w_i)$ because $v \in V(w_i)$, ${\rm deg}(v) = 2$ and
$B_{w_i}$ is primitive. Now, again by Proposition
\ref{generadoresinducido} and Theorem \ref{induce-CI}, for every $u
\in V(G) \setminus \{v, v_1, v_2\}$, we have that $u \in V(w_i)$ if
and only if  $\{B_{w_j}\, \vert \, u \notin V(w_j)\} = \{B_{w_j}\,
\vert \, v,u \notin V(w_j)\}$, or equivalently if $\mu(P_{G
\setminus \{u\}}) = \mu(P_{G \setminus \{u,v\}}) \Leftrightarrow
{\rm ht}(P_{G \setminus \{u\}}) = {\rm ht}(P_{G \setminus
\{u,v\}})$.  Since ${\rm deg}_{G \setminus \{u\}}(v) = 2$, this is
equivalent to $b(G \setminus \{u,v\})
> b(G \setminus \{u\})$.

Finally, if (d) holds, Corollary \ref{2situaciones} and Theorem
\ref{todosgradomayor2} complete the proof. \QED
\end{demo}

\smallskip
This theorem yields Algorithm CI-graph, see Figure \ref{algoritmo},
an algorithm to determine if a graph is a complete intersection.
This method begins by removing all the vertices of degree $1$ and
$2$ iteratively. Whenever we remove a vertex $v$ of degree $2$, we
check whether $b(G) = b(G \setminus \{v\})$. In the positive case,
we construct a set $W \subset V(G)$ and look for an even closed walk
$w$ such that $V(w) = W$. If such a walk does not exist, then $G$ is
not a complete intersection, otherwise we take $w$ a shortest even
walk such that $V(w) = W$ and define the binomial $B_w$; one can
obtain such an even walk in polynomial time by means of the
algorithm proposed in \cite{LaPauPapa}. Once we have removed every
vertex of degree $\leq 2$, either we get a trivial graph or we reach
a graph $G'$ where every vertex has degree $> 2$. If there exists a
connected component of $G'$ which is neither an odd band nor an even
M\"obius band, then $G$ is not a complete intersection. Otherwise we
can construct a set of $r = {\rm ht}(P_G)$ binomials
$\{B_{w_1},\ldots,B_{w_r}\} \subset P_G$, and $G$ is a complete
intersection if and only if $P_G = (B_{w_1},\ldots,B_{w_r})$. For
checking this equality we use Theorem \ref{fs}. It is worth pointing
out that in \cite{F-M-S} the authors give a polynomial algorithm to
decide if a matrix is dominating; thus one can check if the equality
$P_G = (B_{w_1},\ldots,B_{w_r})$ holds in polynomial time.

\medskip As a direct consequence of this algorithm we have the
following result.

\begin{Corollary}The problem of determining whether a graph is a
complete intersection is in the complexity class $\mathcal P$.
\end{Corollary}
\begin{demo}Counting the number of connected components of a graph and deciding
whether a graph is bipartite, and thus computing $b(H)$, can be done
in polynomial time for every graph $H$. To prove the result it only
remains to prove that, given a connected graph $H$, one can decide
if $H$ is either an odd band or an even M\"obius band in polynomial
time. For this purpose we propose a polynomial time algorithm that
returns {\sc True} if $H$ is an odd band or an even M\"obius band,
or {\sc False} otherwise. If $H$ has $4$ vertices then we return
{\sc True} if and only if $H = \mK_4$. If $H$ has more than $4$
vertices, the algorithm lies on the fact if $H$ is an odd band or an
even M\"obius band with $V(H) = \{a_1,\ldots,a_r,b_1,\ldots,b_r\}$
and edges $\{a_i,a_{i+1}\}, \{b_i,b_{i+1}\}$ for all $i \in
\{1,\ldots,r-1\}$, $\{a_i,b_i\}$ for all $i \in \{1,\ldots,r\}$;
then $b(H \setminus \{a_i,b_i\}) = 1$ for all $i \in
\{1,\ldots,r\}$, $b(H \setminus \{a_i,a_{i+1}\}) = 0$ for every $i
\in \{1,\ldots,r-1\}$; moreover $b(H \setminus \{a_1,a_r\}) = 0$ if
$H$ is an odd band, and $b(H \setminus \{a_1,b_r\}) = 0$ if $H$ is
an even M\"obius band. The algorithm receives as input the graph
$H$, if $H$ is not $3$-regular or $H$ is bipartite, we return {\sc
False}. Otherwise we take $a_1 \in V(H)$ an arbitrary vertex and
denote $N_H(a_1) := \{w_1,w_2,w_3\}$. We compute $c_i := b(H
\setminus \{a_1,w_i\})$ for all $i \in \{1,2,3\}$ and assume that
$c_1 \geq c_2 \geq c_3$. If $c_2 \geq 1$ or $c_1 = 0$, then we
return {\sc False}. Otherwise we set $b_1 := w_1$. Now we take $a_2
\in N_H(a_1)$, such that $a_2 \notin H_1 := \{a_1,b_1\}$ and iterate
this process until we get that $V(H) =
\{a_1,\ldots,a_r,b_1,\ldots,b_r\}$, $\{a_i,a_{i+1}\} \in E(H)$ for
all $i \in \{1,\ldots,r-1\}$ and $\{a_i,b_i\} \in E(H)$ for all $i
\in \{1,\ldots,r\}$. Thus, we return {\sc True} if and only if
$\{b_i,b_{i+1}\} \in E(H)$ for all $i \in \{1,\ldots,r-1\}$. \qed
\end{demo}

\begin{figure}

\centering
\begin{tabular}{|p{10cm}|}
\hline
$$\begin{array}{l} \mbox{\bf Algorithm CI-graph } \end{array}$$
$$\begin{array}{cl}
\ \mbox{Input:} & G \mbox{ a simple graph.} \\
\ \mbox{Output:} & \mbox{{\sc True} if } G \mbox{ is a complete intersection or {\sc False} otherwise} \\
\end{array}
$$

\medskip

{
\begin{algorithmic}
\STATE $H := G$; $\mathfrak B := \emptyset$

\WHILE {$\exists\, v \in V(H)$ with ${\rm deg}_{H}(v) \leq 2$}

\IF {${\rm deg}_{H}(v) = 2$ and $b(H \setminus \{v\}) = b(H)$}

\STATE $W := \{v\} \cup N_{H}(v) \cup \{u \in V(H) \  \vert \  b(H
\setminus \{u,v\}) >  b(H \setminus \{u\})\}$

\IF {not exists an even closed walk such that $V(w) = W$}

\RETURN {\sc False}

\ENDIF

\STATE Let $w$ be a shortest even closed walk with $V(w) = W$.

\STATE $\mathfrak B := \mathfrak B \cup \{B_w\}$

\ENDIF

\STATE {$H := H \setminus \{v\}$}

\ENDWHILE

\STATE Let $H_1,\ldots,H_s$ be the connected components of $H$

 \IF {exists $i$ such that $H_i$ is not odd band or even M\"obius band }

\RETURN {\sc False}

\ENDIF

\STATE Let $\mathfrak B_i$ be a minimal set of generators of
$P_{H_i}$ for $1 \leq i \leq s$.

\IF {$P_G = \langle \mathfrak B \cup \mathfrak B_1 \cup \cdots \cup
\mathfrak B_s \rangle$}

\RETURN {\sc True}

\ENDIF

\RETURN {\sc False}
\end{algorithmic}
}

\\

\hline

\end{tabular}
\caption{Pseudo-code for checking whether a graph is a complete
intersection. It returns {\sc True} if $G$ is a complete
intersection and {\sc False} otherwise.} \label{algoritmo}
\end{figure}

\bigskip

Let us illustrate how Algorithm CI-graph works with an example.
\begin{Example}\label{ejemploutil}Let us prove that the graph $G$ in {\rm Figure
\ref{fig5}} is not a complete intersection.
\begin{figure}
\begin{center}
\scalebox{1.2} 
{
\begin{pspicture}(0,-1)(3,1)
{\tiny \psdots[dotsize=0.1](0.15,0.3) \rput(0.15,0.43){$v_2$}
\psdots[dotsize=0.1](0.2,-0.9) \rput(0.2,-1.1){$v_1$}
\psdots[dotsize=0.1](1,-0.3) \rput(1.1,-0.5){$v_3$}
\psdots[dotsize=0.1](1.8,-0.28) \rput(1.75,-0.45){$v_4$}
\psdots[dotsize=0.1](2.6,-0.95) \rput(2.8,-0.95){$v_5$}
\psdots[dotsize=0.1](2.6,0.3) \rput(2.8,0.3){$v_6$}
\psdots[dotsize=0.1](1.4,0.9) \rput(1.4,1){$v_7$}}

\psline[linewidth=0.02cm](0.12,0.3)(1,-0.28)
\psline[linewidth=0.02cm](1,-0.28)(1.80,-0.26)
\psline[linewidth=0.02cm](1.8,-0.26)(2.62,0.3)
\psline[linewidth=0.02cm](2.6,0.3)(2.6,-0.9)
\psline[linewidth=0.02cm](2.6,-1)(1.76,-0.25)
\psline[linewidth=0.02cm](1,-0.28)(0.2,-0.9)
\psline[linewidth=0.02cm](0.2,-0.9)(0.18,0.28)
\psline[linewidth=0.02cm](1,-0.28)(1.4,0.9)
\psline[linewidth=0.02cm](1.42,0.9)(2.6,0.3) {\tiny
\rput(0,-0.155){$e_1$} \rput(0.7,0.1){$e_2$}
\rput(0.75,-0.75){$e_3$} \rput(1.4,-0.15){$e_4$}
\rput(2.1,-0.7){$e_5$} \rput(2.3,-0.1){$e_7$} \rput(2.8,-0.2){$e_6$}
\rput(2,.7){$e_8$} \rput(1.1,.5){$e_9$}}
\end{pspicture}
}
\end{center}
\caption{A non complete intersection graph.} \label{fig5}
\end{figure}
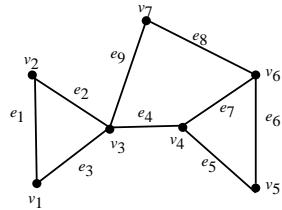

Firstly, we observe that ${\rm deg}_G(v_7) = 2$ and $b(G) = 0 = b(G
\setminus \{v_7\})$. We set $W_1 := \{v_7\} \cup N_G(v_7) \cup \{u
\in V(G) \, \vert\, b(G \setminus \{u,v_7\}) > b(G \setminus
\{u\})\} = \{v_3,v_4,v_6,v_7\}$. We construct $w_1$ a shortest even
closed walk with $V(w_1) = W_1$. Thus $w_1$ is the length four cycle
$w_1 = (v_7, v_3, v_4, v_6, v_7)$ and $B_{w_1} = x_4 x_8 - x_7 x_9$.

Now we consider the graph $H := G \setminus \{v_7\}$ and observe
that ${\rm deg}_H(v_1) = 2$ and $b(H) = 0 = b(H \setminus \{v_1\})$.
We set $W_2 := \{v_1\} \cup N_H(v_1) \cup \{u \in V(H) \, \vert\,
b(H \setminus \{u,v_1\}) > b(H \setminus \{u\})\} = V(H)$. We
construct $w_2$ a shortest even closed walk with $V(w_2) = W_2$.
Then $w_2 = (C_1, \mP_1, C_2, -\mP_1)$ where $C_1 :=
(v_3,v_1,v_2,v_3)$ and $C_2 := (v_4,v_5,v_6,v_4)$ are odd cycles and
$\mP_1$ is the length one path $\mP_1 := (v_3,v_4)$; thus $B_{w_2} =
x_1 x_4^2 x_6 - x_2 x_3 x_5 x_7$.

Hence, we consider the graph $H' := H \setminus \{v_1\}$. We observe
that every $u \in V(H'),\, u \neq v_4$ either has degree $1$ or has
degree $2$ and $b(H') \neq b(H' \setminus \{u\})$. Thus, one can
remove one by one every vertex of $H'$ until getting a trivial
graph.

Then, we have that $G$ is a complete intersection if and only if
$P_G = (B_{w_1}, B_{w_2})$.  We denote by $B$ the $2 \times 9$
matrix $B := \left( \begin{array}{ccccccccc} 0 & 0 & 0 & 1 & 0 & 0 & -1 & 1 & -1 \\
1 & -1 & -1 & 2 & -1 & 1 & -1 & 0 & 0 \end{array} \right)$, then
$\Delta_2(B) = 1$ and it has a square mixed submatrix $B' := \left( \begin{array}{cc} 1 & -1 \\
2 & -1 \end{array} \right)$. Thus $P_G$ is not a complete
intersection.
\end{Example}

\section{Theta graphs and complete intersections} \label{sec5}

This section is devoted to prove that if $G$ is a complete
intersection and there are three paths $\mP_1,\mP_2$ and $\mP_3$ of the
same parity connecting $x, y \in V(G)$ that only meet at their ends,
i.e., $V(\mP_i) \cap V(\mP_j) = \{x,y\}$ for $1 \leq i < j \leq 3$,
then $\mP_1,\mP_2$ and $\mP_3$ are all odd paths and $\{x,y\} \in
E(G)$. We will prove this result in Theorem \ref{notheta}. To prove
this we will first introduce two results concerning the vertices of
degree $2$ in a complete intersection graph, namely Lemma
\ref{2subgrafos} and Proposition \ref{contraccion}. The first one is
a technical lemma which will very useful in the sequel. The second
one describes an operation in $G$ that leads to another graph $G'$
with less vertices and edges than $G$ and $G'$ is a complete
intersection whenever $G$ is.

\begin{Lemma}\label{2subgrafos}Let $v$ be a vertex of degree $2$ and $H_1, H_2$ two
induced subgraphs such that $v \in V(H_i),\, {\rm deg}_{H_i}(v) = 2$
and $b(H_i \setminus \{v\}) = b(H_i)$ for $i = 1,2$. If $G$ is a
complete intersection, then $b(H \setminus \{v\}) = b(H)$, where $H
:= [V(H_1) \cap V(H_2)]$.
\end{Lemma}
\begin{demo}Let $\{B_{w_1},\ldots,B_{w_r}\}$ be a minimal set of generators
of $P_G$ where $w_j$ is an even closed walk for $1 \leq j \leq r$.
For every $G' \in \{G,H_1,H_2,H\}$, by Proposition
\ref{generadoresinducido} and Theorem \ref{induce-CI}, we have that
$G'$ and $G' \setminus \{v\}$ are complete intersections minimally
generated by $\{B_{w_i} \, \vert \, V(w_i) \subset V(G')\}$ and
$\{B_{w_i} \, \vert \, v \notin V(w_i) \subset V(G')\}$,
respectively. Thus, $$1 + b(G') - b(G' \setminus \{v\}) = {\rm
deg}_{G'}(v) - 1 + b(G') - b(G' \setminus \{v\}) = $$  $$= {\rm
ht}(P_{G'}) - {\rm ht}(P_{G' \setminus \{v\}}) = |\{i\, \vert \, v
\in V(w_i) \subset V(G')\}|.$$ In particular, for $i = 1,2$ we have
that there exists a unique $j_i \in \{1,\ldots,r\}$ such that $v \in
V(w_{j_i}) \subset V(H_i)$. We claim that $j_1 = j_2$, indeed $v \in
V(w_{j_1}) \cap V(w_{j_2})$ and $1 \geq 1 + b(G) - b(G \setminus
\{v\}) = |\{i \, \vert \, v \in V(w_i)\}|.$ Thus, $V(w_{j_1})
\subset V(H_1) \cap V(H_2) = V(H)$ and we can conclude that $1 +
b(H) - b(H \setminus \{v\}) = |\{j \, \vert\, v \in V(w_j) \subset
V(H) \}| = 1,$ and $b(H) = b(H \setminus \{v\})$.\QED
\end{demo}

\bigskip

The second result concerns an operation in a graph, which is called
the contraction of a graph in a vertex of degree $2$. We will prove
that if $G$ is a complete intersection, the contraction of $G$ in a
vertex of degree $2$ preserves the property of being a complete
intersection.

\begin{Definition}
Let $G$ be a graph with a vertex $v$ of degree $2$ which does not
belong to a triangle, i.e., $N_G(v) = \{u_1,u_2\}$ and $\{u_1,u_2\}
\notin E(G)$. We define the {\it contraction of $G$ in $v$} as the
graph $G_v^c$ obtained by contracting the two edges incident to $v$.
More precisely, $G^c_v$ is the graph with
$$V(G^c_v) :=  (V(G) \, \setminus \, \{v,u_1,u_2\}) \cup \{u\} \ {\rm and}$$
$$E(G^c_v) := E(G \setminus \{v,u_1,u_2\})  \cup  \ \bigg\{
\{u, x\} \, \mid \, \{u_1,x\} {\rm \ or \ } \{u_2,x\} \in E(G) {\rm
\ and\ } x \neq v\bigg\}.$$
\end{Definition}

\begin{Proposition}\label{contraccion}Let $G$ be a graph with a vertex $v$ of degree
$2$ which does not belong to a triangle. If
 $G$ is a complete intersection, then so is $G_v^c$.
\end{Proposition}
\begin{demo}For every even (respect. odd) closed walk $w = (z_1,\ldots,z_r =
z_1)$ in $G$, we define $\widehat{w}$ as the even (respect. odd)
closed walk in $G_v^c$ constructed as follows. Assume that $z_1
\notin \{v,u_1,u_2\}$, for every $i \in \{2,\ldots,r-1\}$ such that
$z_i = v$, then $z_{i-1},z_{i+1} \in \{u_1,u_2\}$ and we set
$\widehat{w} := (z_1,\ldots,z_{i-2},u,z_{i+2},\ldots,z_r = x_1)$ and
whenever $z_i \in \{u_1,u_2\}$ with $z_{i-1} \neq v$, $z_{i+1} \neq
v$ then we set $\widehat{w} :=
(z_1,\ldots,z_{i-1},u,z_{i+1},\ldots,z_r = z_1)$. Note that it might
happen that $w$ passes by $u_1, u_2$ or $v$ more than once.
Moreover, for every closed walk $w'$ in $G_v^c$ one can find another
$w$ in $G$ such that $w' = \widehat{w}$.

We have that $G$ is bipartite if and only if so is $G_v^c$, indeed
 $V_1, V_2$ is a bipartition for $G$ with
$u_1 \in V_1$ if and only if $V_1', V_2'$ is a bipartition for
$G_v^c$, where $u \in V_1'$, $V_1 \setminus \{u_1,u_2\} = V_1'
\setminus \{u\}$ and $V_2 \setminus \{v\} = V_2'$. Moreover, since
$|V(G_v^c)| = m - 2$ and $|E(G_v^c)| = n - |N_G(u_1) \cap N_G(u_2)|
- 1$, we have that
$$ {\rm ht}(P_{G_v^c}) = {\rm ht}(P_G) - |N_G(u_1) \cap N_G(u_2)| +
1.$$ Moreover $|N_G(u_1) \cap N_G(u_2)| \leq 2$, otherwise $G$ has a
subgraph $\mathcal K_{2,3}$, which is not possible by Corollary
\ref{subgraphK23}. This proof falls naturally into two parts.

If $N_G(u_1) \cap N_G(u_2) = \{v\}$. Assume that $e_{n-1} = \{u_1,
v\},\, e_n = \{u_2, v\}$ and set $e_i' := e_i$ if $e_i \in E(G_v^c)$
and $e_i' := \{u,z\}$ if either $e_i = \{u_1,z\}$ or $e_i =
\{u_2,z\}$ for every $1 \leq i \leq n-2$; then $E(G_v^c) =
\{e_1',\ldots,e_{n-2}'\}$.

Consider now the morphism $\psi: k[x_1,\ldots,x_n] \longrightarrow
k[x_1,\ldots,x_{n-2}]$ induced by $x_{n-1} \mapsto 1$\,, $x_{n}
\mapsto 1$\, and $x_i \mapsto x_i$ for every $i \in
\{1,\ldots,n-2\}$. It is easy to check that for every even closed
walk $w$ in $G$, then $\psi(B_w) = B_{\widehat{w}}$. This implies
that $\psi(P_G) = P_{G_v^c}$. Since ${\rm ht}(P_G) = {\rm
ht}(P_{G_v^c})$, we get that $G_v^c$ is a complete intersection
because if $\mathfrak B$ is a set of generators of $P_G$ then
$\psi(\mathfrak B)$ is a set of generators of $P_{G_v^c}$.

Secondly consider the case where $N_G(u_1) \cap N_G(u_2) = \{v,
z\}$. Suppose that $e_{n-3} = \{u_1, z\}$,  $e_{n-2} = \{u_2, z\}$,
$e_{n-1} = \{u_2, v\}$ and $e_n = \{u_1, v\}$ and set $e_i' := e_i$
if $e_i \in E(G_v^c)$, $e_i' := \{u, t\}$ if either $e_i = \{u_1,
t\}$ or $e_i = \{u_2, t\}$ for all $i \in \{1,\ldots,n-4\}$ and
$e_{n-3}' := \{u,z\}$; then $E(G_v^c) = \{e_1',\ldots,e_{n-3}'\}$.

Consider now the morphism $\psi: k[x_1,\ldots,x_n] \longrightarrow
k[x_1,\ldots,x_{n-3}]$ defined by $x_{n-2} \mapsto x_{n-3}$\,,
$x_{n-1} \mapsto 1$\,, $x_{n} \mapsto 1$\, and $x_i \mapsto x_i$ for
every $i \in \{1,\ldots,n-3\}$. Proceeding as before, we get that
$\psi(P_G) = P_{G_v^c}$.

Suppose that $G$ is a complete intersection and consider the cycle
$w := (v, u_1, z, u_2, v)$, then the quadric $B_w := x_{n-2}  x_{n}
- x_{n-3}  x_{n-1} \in P_G$. Since $P_G$ is a homogeneous ideal
which does not contain any linear form, we get that there exists a
minimal set of generators $\mathfrak B$ of $P_G$ such that
$\mathfrak B = \{B_{w_1},\ldots,B_{w_r}, B_w\}$ for some even closed
walks $w_1,\ldots,w_r$ in $G$ and then $r + 1 = {\rm ht}(P_G)$.
Since $\psi(P_G) = P_{G_v^c}$ and $\psi(B_w) = 0$, it follows that
$\{\psi(B_{w_1}), \psi(B_{w_2}), \ldots, \psi(B_{w_r})\}$ generates
$P_{G_v^c}$ and ${\rm ht}(P_{G_v^c}) = r$, thus it is a complete
intersection. \QED
\end{demo}

\bigskip

\begin{figure}[h]
\begin{center}
\setlength{\unitlength}{.04cm}
\begin{picture}(150,65)(20,10)

\put(50,65){$G$}

\put(12,65){\footnotesize $v$}

\put(12,5){\footnotesize $u_2$}

\put(20,70){\circle*{3}}

\put(20,70){\line(0,-1){60}}

\put(20,70){\line(2,-1){40}}

\put(20,10){\circle*{3}}

\put(20,10){\line(2,1){40}}

\put(40,50){\circle*{3}}

\put(60,53){\footnotesize $u_1$}

\put(40,50){\line(1,0){40}}

\put(60,50){\circle*{3}}

\put(80,50){\circle*{3}}

\put(40,30){\circle*{3}}

\put(40,30){\line(1,0){40}}

\put(40,30){\line(0,1){20}}

\put(80,30){\line(0,1){20}}

\put(60,30){\circle*{3}}

\put(80,30){\circle*{3}}


\put(150,65){$G_v^c$}

\put(160,53){\footnotesize  $u$}

\put(140,50){\circle*{3}}

\put(140,50){\line(1,0){40}}

\put(160,50){\circle*{3}}

\put(160,50){\line(0,-1){20}}

\put(180,50){\circle*{3}}

\put(140,30){\circle*{3}}

\put(140,30){\line(1,0){40}}

\put(140,30){\line(0,1){20}}

\put(180,30){\line(0,1){20}}

\put(160,30){\circle*{3}}

\put(180,30){\circle*{3}}
\end{picture}
\end{center}
\caption{$G$ is bipartite but it is not a ring graph, hence it is
not a complete intersection. Nevertheless, $v$ is a vertex of degree
$2$ which does not belong to a triangle and $G_v^c$ is a complete
intersection.} \label{fig4}
\end{figure}
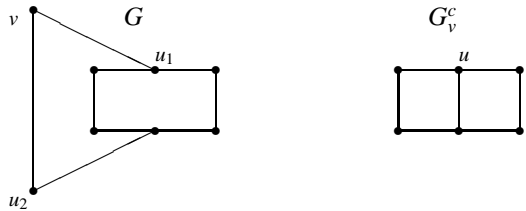

The converse of Proposition \ref{contraccion} is not true in
general, as one can see in Figure \ref{fig4}

Now we introduce the concept of theta graph and use the previous
results to prove Theorem \ref{notheta}, which asserts that odd theta
graph whose base vertices are not adjacent, and also even theta
graphs, are forbidden subgraphs in a complete intersection graph.
This is the main result of this section.

\begin{Definition}\label{thetadefinition}A {\it theta graph} $T$ with {\it base vertices} $x, y$
is a graph with $V(T) := V(\mP_1) \cup V(\mP_2) \cup V(\mP_3)$, where
$\mP_1, \mP_2$ and $\mP_3$ are three paths of length $\geq 2$
connecting $x$ and $y$ such that $V({\mP}_i) \cap V({\mP}_j)=\{x,y\}$
for $1 \leq i < j \leq 3$. If $\mP_1, \mP_2$ and $\mP_3$ are all even
$($respect. odd$)$ paths, then $T$ is called an {\it even}
$($respect. {\it odd}$)$ {\it theta graph}.
\end{Definition}

\begin{Remark}Theta graphs are sometimes defined in the literature
to have $E(T) = E(\mP_1) \cup E(\mP_2) \cup E(\mP_3).$ However, from
our definition we have $E(\mP_1) \cup E(\mP_2) \cup E(\mP_3) \subset
E(T)$; this is, it might have edges connecting a vertex in $\mP_i$
and a vertex in $\mP_j$ with $1 \leq i < j \leq 3$ or even connecting
two vertices $a_i,a_j$ of $\mP_k = (x = a_0, a_1, \ldots, a_{r-1},
a_r = y)$ with $0 \leq i < i+1 < j \leq r$ and $1 \leq k \leq 3$.
\end{Remark}

To prove Theorem \ref{notheta} we need a lemma, whose proof is
almost immediate.

\begin{Lemma}\label{verticeciclopar}Let $v$ be a vertex of degree $2$ and $C$ an even cycle
such that $v \in V(C)$. Then, $G$ is bipartite $\Longleftrightarrow
G \setminus \{v\}$ is bipartite.
\end{Lemma}
\begin{demo}Let $u_1, u_2$ be the neighbors of $v$, then $u_1, u_2 \in V(C)$.
Hence, if $V_1, V_2$ is a bipartition for $G \setminus \{v\}$ and
$u_1 \in V_1$, then necessarily $u_2 \in V_1$, so $V_1, V_2 \cup
\{v\}$ is a bipartition for $G$. The other implication is obvious.
\QED
\end{demo}

\bigskip

\begin{Theorem}\label{notheta}Neither odd theta graphs whose base vertices
are not adjacent nor even theta graphs are complete intersections.
\end{Theorem}
\begin{demo}Suppose that there exists an even theta graph or an odd theta
graph whose base vertices are not adjacent which is a complete
intersection. Let $G$ be the smallest graph with respect to $|V(G)|$
with this property, we denote by $x,y$ the base vertices of $G$ and
$\mP_1, \mP_2, \mP_3$ the three even or odd paths connecting $x$ and
$y$ such that $V(\mP_i) \cap V(\mP_j) = \{x,y\}$ for $1 \leq i < j
\leq 3$.

If $G$ is $3$-regular, then one can write $N_G(x) = \{u_1,u_2,u_3\}$
with $u_i \in V(\mP_i)$ and we claim that $b(G \setminus \{u_i\}) =
b(G \setminus \{x,u_i\})$ for $1 \leq i \leq 3$.  Indeed, we can
suppose that $i = 1$, then $x$ has degree $2$ in $G \setminus
\{u_1\}$ and belongs to the even cycle $(\mP_2, -\mP_3)$; then by
Lemma \ref{verticeciclopar} we have that $b(G \setminus \{u_1\}) =
b(G \setminus \{x,u_1\})$. On the other hand, since $G$ is a
complete intersection $3$-regular graph, it is an odd band or an
even M\"obius band. Thus one can write $V(G) =
\{a_1,\ldots,a_r,b_1,\ldots,b_r\}$ and assume that $x = a_i$ for
some $1 \leq i \leq r$. Then $b(G \setminus \{b_i\}) = 0$ and $b(G
\setminus \{x,b_i\}) = 1$, which is a contradiction.

Then, by Corollary \ref{2situaciones} there exists $v \in V(G)$ of
degree $2$ and we will assume that $v \in V(\mP_3)$. Then we consider
$C_i := (\mP_i, -\mP_3)$ and $H_i := [V(C_i)]$ for $i = 1,2$. Since
$C_i$ is an even cycle and $v \in V(C_i)$, it follows by Lemma
\ref{verticeciclopar} that $b(H_i) = b(H_i \setminus \{v\})$ for $i
= 1,2$. Thus by Lemma \ref{2subgrafos}, $b(H) = b(H \setminus
\{v\})$ where $H := [V(H_1) \cap V(H_2)] = [V(\mP_3)]$, let us prove
that this is not possible. We denote $\mP_3 := (x = a_0, a_1, \ldots,
a_r = y)$, then $v = a_i$ for some $1 \leq i < r$.

If $r = 2$, then $v = a_1$. Moreover, if $\{a_0,a_2\} \in E(G)$ we
have that $b(H) = 0 \neq 1 = b(H \setminus \{v\})$ and if
$\{a_0,a_2\} \notin E(G)$ we have that $b(H) = 1 \neq 2 = b(H
\setminus \{v\})$.

If $r = 3$, then $G$ is an odd theta graph and $\{x, y\} \notin
E(G)$, we can assume that $v = a_1$, then if $\{a_0,a_2\} \in E(G)$
we have that $b(H) = 0 \neq 1 = b(H \setminus \{v\})$ and if
$\{a_0,a_2\} \notin E(G)$ we have that $b(H) = 1 \neq 2 = b(H
\setminus \{v\})$.

If $r \geq 4$, let us prove that $\{a_{i-1},a_{i+1}\} \in E(G)$.
Assume that $\{a_{i-1},a_{i+1}\} \notin E(G)$ and consider $G_v^c$,
which is an even or odd theta graph. By Proposition
\ref{contraccion}, $G_v^c$ is a complete intersection and
$|V(G_v^c)| < |V(G)|$. This contradicts the minimality of $G$ unless
if $G_v^c$ is an odd theta graph whose base vertices are adjacent.
This can only happen if $v = a_1$ and $\{a_2,y\} \in E(G)$, or $v =
a_r$ and $\{a_{r-1},x\} \in E(G)$. Without loss of generality we can
assume that $v = a_1$ and $\{a_2,y\} \in E(G)$. Then, the induced
subgraph $G' := [V(\mP_1) \cup V(\mP_2) \cup V(\mP_3')]$ with $\mP_3' =
(x, a_1, a_2, y)$ is an odd theta graph whose base vertices $x,y$
are not adjacent and is a complete intersection, but this is not
possible by the minimality of $G$. Thus, $\{a_{i-1},a_{i+1}\} \in
E(G)$.

Since $\{a_{i-1},a_{i+1}\} \in E(G)$, then $b(H) = 0$; let us see
that $H \setminus \{v\}$ is bipartite and $b(H \setminus \{v\}) =
1$. Indeed, suppose that $H \setminus \{v\}$ is not bipartite, we
denote $a_j' := a_j$ for $1 \leq j \leq i-1$, $a_j' := a_{j+1}$ for
$i \leq j < r$. Since $\{a_j',a_{j+1}'\} \in E(H \setminus \{v\})$
for $1 \leq j < r$, there exists $1 \leq j < k < r$ such that
$\{a_j',a_k'\} \in E(H)$ and $j \equiv k\ ({\rm mod}\, 2)$. We
separate three cases:
\begin{itemize}[leftmargin=.8cm] \item[(a)] if $k < i$, then
$\{a_j, a_k\} \in E(H)$.
\item[(b)] if $j \geq i$, then
 $\{a_{j+1}, a_{k+1}\} \in E(H)$
\item[(c)] $j < i \leq k$, then $\{a_j, a_{k+1}\} \in
E(H)$.
\end{itemize}
If (a) holds, we denote $\mP_3' = (x = a_0,\ldots, a_j, a_k,\ldots,
a_{i-1}, a_{i+1},\ldots, a_r = y)$. If (b) holds, we denote $\mP_3' =
(x = a_0,\ldots, a_{i-1}, a_{i+1},\ldots, a_{j+1}, a_{k+1},\ldots,
a_r = y)$. If (c) holds, we denote $\mP_3' = (x = a_0,\ldots, a_j,
a_{k+1},\ldots, a_r = y)$. In the three cases $\mP_3'$ is an walk of
the same parity of $\mP_3$ connecting $x$ and $y$ with $V(\mP_3')
\subsetneq V(\mP_3)$, but this contradicts again the minimality of
$G$.\QED
\end{demo}

\section{Structure theorems for complete intersection graphs} \label{sec7}

The goal of this section is to prove two structure theorems for
complete intersection connected graphs; namely Theorem
\ref{teoremaestructura} and Theorem \ref{ultimo}. Given a graph, it
can be partitioned into
 two induced subgraphs $C$ and $R$, such that $V(C) = V(C_1)
\bigsqcup \cdots \bigsqcup V(C_s)$ where $C_1,\ldots,C_s$ are odd
primitive cycles and $R$ is a bipartite graph. Note that this
partition might not be unique and when $G$ is bipartite, one has
that $C$ is the empty graph. Whenever we have a partition with these
properties we write $G = [C; R]$. In order to characterize the
complete intersection property on $G$, we propose to characterize
when $C$ and $R$ are complete intersections, and then determine the
admissible edges connecting $C$ and $R$. Since $R$ is a bipartite
graph it turns out that it is a complete intersection if and only if
it is a ring graph (see \cite[Corollary 3.3]{Ring}). Theorem
\ref{teoremaestructura} will give necessary conditions for a
connected graph to be a complete intersection by determining when
$C$ is a complete intersection. Finally, if $C$ is connected and $R$
is $2$-connected, Theorem \ref{ultimo} characterizes the complete
intersection property by obtaining all possible edges connecting $C$
and $R$.

\smallskip Let us start with this proposition.

\begin{Proposition}\label{3impares}Let $G$ be a complete intersection connected graph,
then there are at most two vertex disjoint odd cycles in $G$.
\end{Proposition}

We need the following technical result which is included in the
proof of \cite[Theorem 5.3]{Tatakis-Thoma}. Recall that a {\it
block} is a maximal connected subgraph $B$ of $G$ such that if one
removes any of its vertices it is still connected. A graph is {\it
$2$-connected} if it only has one block and more than $2$ vertices.

\begin{Lemma}\label{tt} Let $G$ be a complete intersection
$2$-connected graph and let $C_1, C_2$ be two odd cycles in
$G$.\begin{itemize}[leftmargin=1cm] \item[{\rm (a)}] If $V(C_1)
\cap  V(C_2) = \{v\}$, then there exists an $e \in E(G)$ such that
$v \notin e$ and $e \cap V(C_i) \not= \emptyset$ for $i \in
\{1,2\}$.
\item[\rm (b)] If $C_1$ and $C_2$ are vertex disjoint, then there exist
$e_1, e_2 \in E(G)$ such that $e_1 \cap e_2 = \emptyset$ and $e_i
\cap V(C_j) \not= \emptyset$ for $i,j \in \{1,2\}$.
\end{itemize}
\end{Lemma}

\demos Assume that $G$ is a complete intersection with three vertex
disjoint odd cycles and let $G'$ be the smallest connected induced
subgraph with this property. We denote by $C_1 = (a_1,\ldots,
a_{r_1}, a_1)$, $C_2 = (b_1,\ldots, b_{r_2}, b_1)$ and $C_3 =
(c_1,\ldots, c_{r_3}, c_1)$ three vertex disjoint odd primitive
cycles of $G'$. By \cite[Theorem 4.2]{Tatakis-Thoma}, $G'$ has
either one or two non bipartite blocks and the proof falls naturally
in two cases.

If $G'$ has only one non bipartite block, then $C_1,  C_2$ and $C_3$
belong to it and, by Lemma \ref{tt}, for every $1 \leq i < j \leq 3$
there exist two edges connecting a vertex of $C_i$ and a vertex of
$C_j$; thus $G' = [V(C_1) \cup V(C_2) \cup V(C_3)]$. $G'$ can not be
a band or a M\"obius band because there are three vertex disjoint
odd primitive cycles, then there exists a vertex $z \in V(G')$ of
degree $2$. Suppose that $z \in V(C_3)$ and denote by $H_i :=
[V(C_i) \cup V(C_3)]$ for $i = 1,2$. We have that $H_i \setminus
\{z\}$ is connected and $V(C_i) \subset V(H_i \setminus \{z\})$,
then $b(H_i) = b(H_i
 \setminus \{z\}) = 0$ for $i = 1,2$. By Lemma \ref{2subgrafos} we
have that $b(H) = b(H \setminus \{z\})$ where $H := [V(H_1) \cap
V(H_2)] = [V(C_3)]$. However, $C_3$ is an odd primitive cycle, then
$b(H) = 0$ and $b(H \setminus \{z\}) = 1$, which is a contradiction.

If $G'$ has two non bipartite blocks, then two of the odd cycles
belong to the same block of $G'$, say $C_1$ and $C_2$. By Lemma
\ref{tt}, $C_1$ and $C_2$ are connected by at least two edges.
Moreover, $C_3$ is not in the same block of $C_1$ and $C_2$. Then we
set $G_1 := [V(C_1) \cup V(C_2)]$ and take $\mP$ a path in $G'$ of
minimum length connecting a vertex of $C_3$ and a vertex of $G_1$.
By the minimality of $G'$ we have that $G' = [V(G_1) \cup V(\mP) \cup
V(C_3)]$. Moreover, we can assume that there exists $s \geq 0$ such
that $\mP = (c_1 = u_0,  u_1,\ldots, u_s, a_1)$; since $\mP$ has
minimum length one can deduce that $u_i \notin V(G_1) \cup V(C_3)$
for $1 \leq i \leq s$, $\{c_j,u_i\} \notin E(G')$ for every $1 \leq
j \leq r_3$, $i > 1$ and $\{a_j,u_i\} \notin E(G')$ for every $1
\leq j \leq r_1$, $i < s$.

Firstly assume that ${\rm deg}_{G'}(c_j) = 2$ for every $j > 1$ and
take $u := c_2$. We set $\mP'$  the shortest path in $G'$ connecting
$c_1$ with a vertex of $C_2$. Then, we can assume that $\mP' = (c_1 =
v_0, v_1,\ldots, v_t = b_1)$ and have that
\begin{itemize}[leftmargin=.8cm]
\item[(a)] $v_i = u_i$ for $0 \leq i \leq s$, $v_i \in V(C_1)$ for $s < i \leq t-1$
and $v_t \in V(C_2)$,
\item[(b)] $\{v_i, v_j\} \notin E(G')$ for $0 \leq i < i+1 < j \leq
t$.
\end{itemize}

 Now we set $H_1 := [V(C_1) \cup V(\mP)
\cup V(C_3)] $ and $H_2 := [V(C_2) \cup V(\mP') \cup V(C_3)]$;
clearly $b(H_i) = b(H_i \setminus \{u\}) = 0$ for $i = 1,2$.
However, if we set $H := [V(H_1) \cap V(H_2)] = [V(C_3) \cup
\{v_1,\ldots,v_{t-1}\}]$, then $V(C_3) \subset V(H)$, ${\rm
deg}_H(c_1) = 3,\, {\rm deg}_H(v_{t-1}) = 1$ and ${\rm deg}_H(v) =
2$ for the rest of vertices of $H$. Thus $b(H) = 0$ and $b(H
\setminus \{u\}) = 1$, which contradicts Lemma \ref{2subgrafos}.

So assume that ${\rm deg}_{G'}(c_j) > 2$ for some $j > 1$ and let us
see that $s = 0$; i.e., $\{a_1,c_1\} \in E(G')$. Indeed, if $s \geq
1$, by the minimality of $G'$ we have that $\{u_1,c_j\} \in E(G')$,
but then there exists an odd primitive cycle $C'$ such that $u_1 \in
V(C') \subset V(C_3) \cup \{u_1\}$, which contradicts the minimality
of $G'$. Thus, $s = 0$, $\{c_1,a_1\} \in E(G')$ and if $\{c_k,a_i\}
\in E(G')$ then $i = 1$ because $G'$ has two blocks. Let $C_3'$ be
an odd primitive cycle such that $a_1 \in V(C_3') \subset V(C_3)
\cup \{a_1\}$ and take $u \in V(C_3'),\, u \neq a_1$. We set $\mP'$
the shortest path in $G'$ connecting $a_1$ with a vertex of $C_2$.
Then, we can assume that $\mP' = (a_1 = v_0, v_1,\ldots,v_t = b_1)$
and we have that
\begin{itemize}[leftmargin=.8cm]
\item[(a)] $v_i \in V(C_1)$ for $0 \leq i \leq t-1$
\item[(b)] $\{v_i, v_j\} \notin E(G')$ for $0 \leq i < i+1 < j \leq
t$.
\end{itemize}

Now we set $H_1 := [V(C_1) \cup V(C_3')] $ and $H_2 := [V(C_2) \cup
V(\mP') \cup V(C_3')]$; and have that $u$ has degree $2$ in $[V(H_1)
\cup V(H_2)]$ and clearly ${\rm deg}_{H_i}(u) = 2$, $b(H_i) = b(H_i
\setminus \{u\}) = 0$ and ${\rm deg}_{H_i}(u) = 2$ for $i = 1,2$.
However, if we set $H := [V(H_1) \cap V(H_2)] = [V(C_3') \cup
\{v_1,\ldots,v_{t-1}\}]$, then $V(C_3') \subset V(H)$, ${\rm
deg}_H(a_1) = 3,\, {\rm deg}_H(v_{t-1}) = 1$ and ${\rm deg}_H(v) =
2$ for the rest of vertices of $H$. Thus $b(H) = 0$ and $b(H
\setminus \{u\}) = 1$, which again contradicts Lemma
\ref{2subgrafos}. \QED

\bigskip

Now that we know that there are at most two vertex disjoint odd
primitive cycles in a complete intersection connected graph, let us
determine how two such cycles can be connected.

\begin{Definition} A graph $G$ is called an {\it odd partial
band} if there exist two vertex disjoint odd primitive cycles
$C_1=(a_1,\ldots,a_{r_1}, a_1)$ and $C_2=(b_1,\ldots, b_{r_2}, b_1)$
such that $V(G) = V(C_1) \cup V(C_2)$ and $$E(G) = E(C_1) \cup
E(C_2) \cup \{\{a_{j_1},b_{k_1}\},\ldots,\{a_{j_s},b_{k_s}\}\},$$
where $s \geq 1,\ 1 \leq j_1 \leq \cdots \leq j_s \leq r_1$, $1 \leq
k_1 \leq \cdots \leq k_s \leq r_2$ and $j_i \equiv k_i\ ({\rm mod}\
2)$ (see {\rm Figure \ref{opb}}).

\begin{figure}
\begin{center}
\scalebox{.7} {
\begin{pspicture}(2,0)(6.5,4)
\psarc[linewidth=0.04](4.27,1.7){2.4}{20}{160}
\psline[linewidth=0.04cm](1.96,2.44)(6.5,2.46)

\psarc[linewidth=0.04](4.34,2.64){2.8}{214.1597}{324.61972}
\psline[linewidth=0.04cm](2.02,1.04)(6.62,1.04)
\psdots[dotsize=0.16](2.04,1.04) \psdots[dotsize=0.16](6.58,1.04)
\psdots[dotsize=0.16](3.18,1.06) \psdots[dotsize=0.16](4.26,1.02)
\psdots[dotsize=0.16](5.38,1.02) \psdots[dotsize=0.16](1.94,2.44)
\psdots[dotsize=0.16](2.66,2.44) \psdots[dotsize=0.16](3.48,2.44)
\psdots[dotsize=0.16](4.28,2.44) \psdots[dotsize=0.16](6.48,2.48)
\psdots[dotsize=0.16](5.06,2.44) \psdots[dotsize=0.16](5.78,2.46)
\psline[linewidth=0.04cm](1.92,2.42)(2.04,1.04)
\psline[linewidth=0.04cm](2.04,1.04)(3.48,2.44)
\psline[linewidth=0.04cm](3.2,1.08)(5.78,2.48)
\psline[linewidth=0.04cm](5.78,2.48)(5.38,1.06)
\psline[linewidth=0.04cm](6.48,2.5)(6.56,1.02)

\rput(1.75,2.65){$a_1$} \rput(2.9,2.65){$a_2$}

\rput(3.5876563,2.65){$a_3$} \rput(4.3009377,2.65){$a_4$}
\rput(5.0895314,2.65){$a_5$}

\rput(5.6,2.65){$a_6$} \rput(6.75,2.65){$a_7$}
\rput(1.966875,0.81){$b_1$} \rput(3.1585937,0.79){$b_2$}
\rput(4.2876563,0.75){$b_3$} \rput(5.4409375,0.73){$b_4$}
\rput(6.7095313,0.71){$b_5$}
\end{pspicture}
} \caption{An odd partial band}\label{opb}
\end{center}
\end{figure}
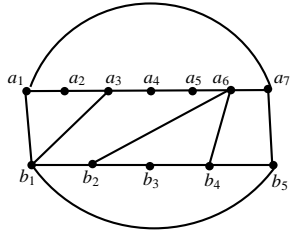
\end{Definition}

\medskip

\begin{Proposition}\label{oddpartilband}Let $G$ be a connected
graph with $V(G) = V(C_1) \cup V(C_2)$ where $C_1$ and $C_2$ are two
vertex disjoint odd primitive cycles. Then, $G$ is a complete
intersection if and only if $G$ is an odd partial band.
\end{Proposition}
\begin{demo}$(\Rightarrow)$ We will proceed by induction on $|V(G)|$.
If $|V(G)| = 6$, i.e., $C_1$ and $C_2$ are triangles, the number of
edges has to be less or equal to $9$, because $2 |E(G)| \leq 3
|V(G)| = 18$ by Corollary \ref{cotasuperior}. If $7 \leq |E(G)| \leq
8$, then one see at once that $G$ is always an odd partial band. If
$|E(G)| = 9$, then $G$ is an odd partial band unless if one can
write $C_1 := (a_1, a_2, a_3, a_1)$,\, $C_2 := (b_1,b_2,b_3,b_1)$
and $E(G) = E(C_1) \cup E(C_2) \cup \{e_1,e_2,e_3\}$ where $e_1 =
\{a_1,b_1\}, e_2 = \{a_1,b_3\}$ and
\begin{enumerate}
\item[(a)] $e_3 = \{a_1,b_2\}$, or
\item[(b)] $e_3 = \{a_2,b_2\}$.
\end{enumerate}
If (a) holds, then ${\rm deg}(a_2) = 2$ and setting $H_i :=
[\{a_1,a_2,a_3,b_i,b_3\}]$ we observe that $b(H_i) = b(H_i \setminus
\{a_2\})$ for $i = 1,2$. However, setting $H := [V(H_1) \cap V(H_2)]
= [\{a_1,a_2,a_3,b_3\}]$, one gets that $b(H) = 0$ and $b(H
\setminus \{a_2\}) = 1$; which contradicts Lemma \ref{2subgrafos}.
If (b) holds, then $G$ has a subgraph $\mK_{2,3}$ with vertices
$\{a_1,a_2,b_1,b_2,b_3\}$, a contradiction to Corollary
\ref{subgraphK23}.

Assume now that $|V(G)| > 6$. Note that $G$ cannot be an even
M\"obius band because $G$ has two vertex disjoint odd primitive
cycles. Thus, if ${\rm deg}(v) > 2$ for every $v \in V(G)$, then $G$
is an odd band, which in particular is an odd partial band.

If every $x \in V(G)$ with ${\rm deg}(x) = 2$ belongs to a triangle,
then we can assume that ${\rm deg}(y) \geq 3$ for every $y \in
V(C_1)$, $C_1$ is not a triangle and $C_2$ is a triangle. Then
$|V(C_1)| \geq 5$, so $2\, |E(G)| = \sum_{v \in V(G)} {\rm deg}(v) =
\sum_{v \in V(C_1)} {\rm deg}(v) + \sum_{v \in V(C_2)} {\rm deg}(v)
\geq 3\, |V(C_1)| + 2\, |V(C_2)| + 5
> 3 |V(G)|$, a contradiction to Corollary \ref{cotasuperior}.

We can suppose that there exists $x \in V(G)$ of degree $2$ which
does not belong to a triangle, we assume that $x \in V(C_1)$. Then,
we consider $G_x^c$ and we have that $V(G_x^c) =
\{a_1',\ldots,a_{r_1-2}',b_1,\ldots,b_{r_2}\}$, where $C_1' :=
(a_1',\ldots,a_{r_1-2}',a_1')$ and $C_2 := (b_1,\ldots,b_{r_2},b_1)$
are odd primitive cycles. Then there exists $2 \leq i \leq r_1 - 1$
such that $x = a_i$,\, $C_1 = (a_1,\ldots,a_{r_1},a_1)$ and
\begin{center} $\{a_j',b_k\} \in E(G_x^c) \Longleftrightarrow \left\{ \begin{array}{cccl} \{a_j,b_k\} \in
E(G) &$  and $& j \leq i-1, &$ or $\\ \{a_{j+2},b_k\} \in E(G) &$
and $& j \geq i-1.  \end{array} \right.$ \end{center} By induction
hypothesis $G_x^c$ is an odd partial band. Thus for every
$\{a_i,b_j\} \in E(G)$ then $i \equiv j\, ({\rm mod}\, 2)$, and for
every $\{a_{i_1}, b_{j_1}\},\, \{a_{i_2}, b_{j_2}\} \in E(G)$ such
that $i_1 \neq i-1$ or $i_2 \neq i+1$, if $i_1 < i_2$, then $j_1
\leq j_2$. Then, $G$ is an odd partial band unless if there exist
two edges $\{a_{i-1},b_{j_1}\},\, \{a_{i+1},b_{j_2}\} \in E(G)$ with
$j_1 > j_2$.

If there exist two adjacent vertices in $C_1$ of degree $2$, then we
take $x = a_i$ one of these two vertices and there can not exist
$\{a_{i-1},b_{j_1}\},\, \{a_{i+1},b_{j_2}\} \in E(G)$ because either
${\rm deg}(a_{i-1}) = 2$ or ${\rm deg}(a_{i+1}) = 2$; hence $G$ is
an odd partial band.

If we are not in the previous situation, then there are at least
three vertices of degree $\geq 3$ in $C_1$. Take $x = a_i \in
V(C_1)$ a vertex of degree $2$ and assume that there exist two edges
$\{a_{i-1},b_{j_1}\},\, \{a_{i+1},b_{j_2}\} \in E(G)$ with $j_1 >
j_2$, let us prove that $G$ is not a complete intersection. Set $u
:= b_{j_2+1}$, we claim that ${\rm deg}(u) = 2$. Indeed $\{a_{i-1},
u\},\,\{a_{i+1}, u\} \not\in E(G)$ because $j_2 + 1 \not\equiv j_2
\equiv i-1 \equiv i+1 \, ({\rm mod}\, 2)$, $\{a_k, u\} \notin E(G)$
if $k < i-1$ because $\{a_{i+1},b_{j_2}\} \in E(G)$ and $j_2 < j_2 +
1$; furthermore $\{a_k, u\} \notin E(G)$ if $k > i+1$ because
$\{a_{i-1},b_{j_1}\} \in E(G)$ and $j_2 + 1 < j_1$.

Take $i' \notin \{i-1, i+1\}$ such that ${\rm deg}(a_{i'}) \geq 3$,
we will assume that $i' < i-1$ and set $j' := {\rm max}\{j \, \vert
\, \{a_{i'},b_j\} \in E(G)\}$. Then necessarily $i' \equiv j' \,
({\rm mod}\, 2)$ and $j' \leq j_2$. Now we consider the even cycle
 $$w_1 :=
(b_{j'}, b_{j'+1},\ldots,b_{j_1}, a_{i-1},
a_{i-2},\ldots,a_{i'},b_{j'})$$ and the even closed walk
$$w_2 :=
(b_{j'}, b_{j'-1},\ldots,b_1, b_{r_2},\ldots,b_{j_2}, a_{i+1},
a_{i+2},\ldots, a_{r_1}, a_1,\ldots,a_{i'},  b_{j'}),$$ which
consists on an even cycle if $j' \neq j_2$ or two odd cycles with
the vertex $b_{j_2}$ in common if $j' = j_2$. For $i = 1,2$ we
denote $H_i := [V(w_i)]$, clearly ${\rm deg}_{H_i}(u) = 2$ and
$b(H_i) = b(H_i \setminus \{u\})$. Let us prove that $b(H) \neq b(H
\setminus \{u\})$ where $H := [V(H_1) \cap V(H_2)] = [ \{a_{i'},
b_{j'}, b_{j_2},\ldots, b_{j_1}\}]$. Indeed, if $j' < j_2 - 1$, then
the vertices $a_{i'}, b_{j'}$ joined by an edge form a connected
component of $H$ , thus $b(H) = 2$ and $b(H \setminus \{v\}) = 3$,
and if $j' = j_2$ or $j' = j_2 - 1$, then $b(H) = 1$ and $b(H
\setminus \{v\}) = 2$. In both cases this is a contradiction to
Lemma \ref{2subgrafos}.

$(\Leftarrow)$ If $G$ is an odd partial band, then ${\rm ht}(P_G) =
s$. We set $e_i := \{a_{j_i}, b_{k_i}\}$ for $1 \leq i \leq s$,
$e_{s+i} := \{a_i, a_{i+1}\}$ for $1 \leq i < r_1$, $e_{s + r_1} :=
\{a_1, a_{r_1}\}$, $e_{s + r_1 + i} := \{b_i, b_{i+1}\}$ for $1 \leq
i < r_2$ and $e_{s + r_1 + r_2} := \{b_1, b_{r_2}\}$.

For every $i \in \{1,\ldots,s-1\}$, let $w_i$ be the even primitive
cycle
$$w_i := (b_{k_i},a_{j_i},a_{j_i +
1},\ldots,a_{j_{i+1}},b_{k_{i+1}}, b_{k_{i+1} - 1},\ldots,
b_{k_{i}})$$ and $w_s := (b_{k_s}, a_{j_s}, a_{j_s + 1},\ldots,
a_{r_1}, a_1, \ldots, a_{j_1}, b_{k_1},\ldots, b_1, b_{r_2},\ldots,
b_{k_s})$.

For every $i \in \{1,\ldots,s\}$ if we write $B_{w_i} = x^{\alpha_i}
- x^{\beta_i}$ with $\alpha_i, \beta_i \in \N^{s+ r_1 + r_2}$, then
$x_i \mid x^{\alpha_i}$. Moreover, for every $1 \leq i < s$,
$x_{i+1} \mid x^{\alpha_i}$ if $j_{i+1} - j_i$ is odd and $x_{i+1}
\mid x^{\beta_i}$ otherwise, and $x_1 \mid x^{\alpha_s}$ if $r_1 -
j_s + j_1$ is odd and $x_1 \mid x^{\beta_s}$ otherwise.

We denote by $B$ the $s \times (s + r_1 + r_2)$ matrix whose $i$-th
row is $\alpha_i - \beta_i$. Note that $\Delta_s(B) = 1$ and for
every $j > s$, the $j$-th column of $B$ has only one nonzero entry
that can be either $+1$ or $-1$. Hence, by Lemma \ref{dominating},
$B$ is dominating if and only if $B'$ is dominating, where $B'$ is
the $s \times s$ submatrix consisting on the first $s$ columns of
$B$. If we denote $B' := (b_{i,j})_{1 \leq i, j \leq s}$, then
$b_{i,j} \neq 0$ if and only if $j - i \in \{0, 1\}$ or $i = s$ and
$j = 1$; thus if there exists a square submatrix $C$ of $B'$ which
is mixed, then $C = B'$ necessarily. Let us see that $B'$ is not
mixed; indeed, $r_1$ is odd and, if we denote by $l_1 = j_2 - j_1,\,
l_2 = j_3 - j_2,\ldots,\,l_{s-1} = j_s - j_{s-1}$ and $l_s = r_1 -
j_s + j_1$; then $r_1 = l_1 + \cdots + l_s$; thus there exists $1
\leq i \leq s$ such that $l_i$ is odd. Then the two nonzero entries
in the $i$-th row of $B'$ are positive and $B'$ is not mixed.
Following Theorem \ref{fs} we can conclude that $P_G$ is a complete
intersection minimally generated by $\{B_{w_1},\ldots,B_{w_s}\}$.
\QED
\end{demo}

\smallskip Now we can state the following structure theorem for
complete intersection graphs.

\begin{Theorem}\label{teoremaestructura}Let $G = [C; R]$ be a complete intersection connected
graph. Then,
\begin{itemize}
\item $R$ is a ring graph
\item $C$ is either the empty graph, an odd primitive cycle, an odd partial band or $C$ has two connected
components which are odd primitive cycles.
\end{itemize}
\end{Theorem}
\begin{demo}This result is a consequence of the fact that the complete
intersection property is hereditary (Theorem \ref{induce-CI}), which
allows us to claim that if $G = [C;R]$ is a complete intersection
graph, then both $R$ and $C$ are complete intersection graphs. Thus,
by Corollary 3.3 in \cite{Ring} it follows that $R$ is a ring graph
and by Proposition \ref{3impares} and Proposition
\ref{oddpartilband} it follows that $C$ is either the empty graph,
an odd primitive cycle, an odd partial band or $C$ has two connected
components which are odd primitive cycles. \QED
\end{demo}

\smallskip

The converse of this statement is not true in general, as the graph
in Figure \ref{figcontraej} shows.

\begin{Example}\label{contraej} Let $G = [C;R]$ be the graph in {\rm Figure
\ref{figcontraej}}, where $C = [\{v_1,v_2,v_3\}]$ is an odd
primitive cycle and $R = [\{v_4,v_5,v_6,v_7\}]$ is a bipartite ring
graph. $R$ is a ring graph and $C$ is a complete intersection,
nevertheless $G$ is not a complete intersection because it contains
$[V(R) \cup \{v_3\}]$, which is $\mK_{2,3}$, as a subgraph; see
Corollary \ref{subgraphK23}.

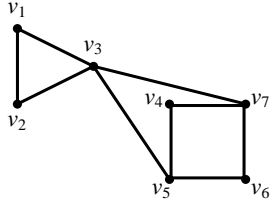
\begin{figure}
\begin{center}
\scalebox{1} {
\begin{pspicture}(2,-0.6)(5,2)
\psframe[linewidth=0.04,dimen=outer](4,-0.25)(5,0.75)
\psdots[dotsize=0.12](4,0.75) \psdots[dotsize=0.12](5,0.75)
\psdots[dotsize=0.12](4,-0.25) \psdots[dotsize=0.12](5,-0.25)
\psdots[dotsize=0.12](2,1.75) \psdots[dotsize=0.12](3,1.25)
\psdots[dotsize=0.12](2,0.75)

\psline[linewidth=0.04cm](2,0.75)(2,1.75)
\psline[linewidth=0.04cm](2,0.75)(3,1.25)
\psline[linewidth=0.04cm](2,1.75)(3,1.25)
\psline[linewidth=0.04cm](3,1.25)(4,-0.25)
\psline[linewidth=0.04cm](3,1.25)(5,0.75)

\rput(2,2){\small $v_1$} \rput(2,0.5){\small $v_2$}
\rput(3,1.5){\small $v_3$} \rput(3.8,0.8){\small $v_4$}
\rput(3.9,-0.4){\small $v_5$} \rput(5.2,-0.4){\small $v_6$}
\rput(5.2,0.8){\small $v_7$}
\end{pspicture}
} \caption{Non complete intersection graph satisfying the conditions
in Theorem \ref{teoremaestructura}.}\label{figcontraej}
\end{center}
\end{figure}
\end{Example}

Under the hypotheses that $R$ is $2$-connected and $C$ is connected,
one has the characterization given in Proposition
\ref{propestructura}. In particular, this proposition states that
there are either $1$ or $2$ vertices in $R$ such that every edge
connecting $R$ and $C$ is incident to one of these vertices. To
state Proposition \ref{propestructura} we need a definition.

\begin{Definition}$G$ is a {\it $1$-clique-sum} of two graphs $G_1$ and $G_2$ if it is obtained by
identifying a vertex $v_1$ of $G_1$ and a vertex $v_2$ of $G_2$.
Analogously, a {\it $2$-clique-sum} of $G_1$ and $G_2$ is obtained
by identifying an edge $e_1$ of $G_1$ and an edge $e_2$ of $G_2$.

\begin{figure}
\begin{center}
\scalebox{1.1} {
\begin{pspicture}(-4,-0.6)(13.762813,1.16375)
\psframe[linewidth=0.04,dimen=outer](2,1)(1,0)
\psdots[dotsize=0.12](1,0) \psdots[dotsize=0.12](2,1)
\psdots[dotsize=0.12](2,0) \psdots[dotsize=0.12](1,1)

\psdots[dotsize=0.12](-1,0) \psdots[dotsize=0.12](-1,1)
\psdots[dotsize=0.12](0,0.5)

\psline[linewidth=0.04cm](-1,0)(-1,1)
\psline[linewidth=0.04cm](-1,0)(0,0.5)
\psline[linewidth=0.04cm](0,0.5)(-1,1)

\psframe[linewidth=0.04,dimen=outer](4,-0.25)(5,0.75)
\psdots[dotsize=0.12](4,0.75) \psdots[dotsize=0.12](5,0.75)
\psdots[dotsize=0.12](4,-0.25) \psdots[dotsize=0.12](5,-0.25)
\psdots[dotsize=0.12](3,0.25) \psdots[dotsize=0.12](3,1.25)

\psline[linewidth=0.04cm](3,0.25)(3,1.25)
\psline[linewidth=0.04cm](3,0.25)(4,0.75)
\psline[linewidth=0.04cm](3,1.25)(4,0.75)

\psframe[linewidth=0.04,dimen=outer](7,1)(8,0)
\psdots[dotsize=0.12](7,1) \psdots[dotsize=0.12](7,0)
\psdots[dotsize=0.12](8,0) \psdots[dotsize=0.12](8,1)
\psdots[dotsize=0.12](6,0.5) \psline[linewidth=0.04cm](7,1)(6,0.5)
\psline[linewidth=0.04cm](7,1)(7,0)
\psline[linewidth=0.04cm](7,0)(6,0.5) \rput(-0.5,-0.3){\footnotesize
$G_1$} \rput(1.5,-0.3){\footnotesize $G_2$}
\rput(4,-.5){\footnotesize $G_3$} \rput(7.2,-0.5){\footnotesize
$G_4$}
\end{pspicture}
} \caption{The graph $G_3$ is a $1$-clique-sum of $G_1$ and $G_2$,
whereas $G_4$ is a $2$-clique-sum of $G_1$ and $G_2$}\end{center}
\end{figure}
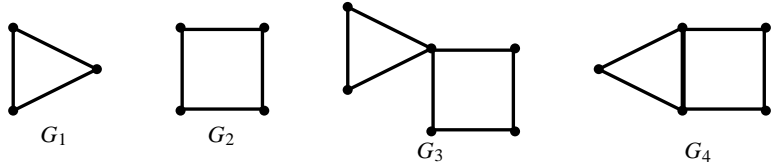
\end{Definition}

\begin{Proposition}\label{propestructura}Let $G = [C;R]$ be a connected graph such that
$R$ is $2$-connected and $C$ is connected. Then, $G$ is a complete
intersection if and only if $R$ is a ring graph and either there
exists $u_1 \in V(R)$ such that $G$ is a $1$-clique-sum of $R$ and
$[V(C) \cup \{u_1\}]$, with $[V(C) \cup \{u_1\}]$ a complete
intersection, or there exist two adjacent vertices $u_1, u_2 \in
V(R)$ such that $G$ is a $2$-clique-sum of $R$ and $[V(C) \cup
\{u_1,u_2\}]$, with $[V(C) \cup \{u_1,u_2\}]$ a complete
intersection.
\end{Proposition}

This result is an immediate consequence of Lemma
\ref{2verticesvecinos} and Lemma \ref{12sumas}.

\begin{Lemma}\label{2verticesvecinos}Let $G = [C; R]$
connected graph such that $R$ is $2$-connected and $C$ is connected.
If $G$ is a complete intersection, then either there exists a vertex
$u_1 \in V(R)$ such that $G$ is a $1$-clique-sum of $R$ and $[V(C)
\cup \{u_1\}]$, or there exist two adjacent vertices $u_1,u_2 \in
V(R)$ such that $G$ is a $2$-clique-sum of $R$ and $[V(C) \cup
\{u_1,u_2\}]$.
\end{Lemma}
\begin{demo}Assume that $G = [C;R]$ is a complete intersection where $R$
is $2$-connected and $C$ is connected. By Theorem
\ref{teoremaestructura}, $C$ is either the empty graph, an odd
primitive cycle or an odd partial band and $R$ is a $2$-connected
ring graph.

Suppose that there exist two edges $e_1 = \{u_1, v_1\}$, $e_2 =
\{u_2, v_2\}$ such that $u_1 \neq u_2,\, u_1, u_2 \in V(R)$, $v_1,
v_2 \in V(C)$ and $u_1$ and $u_2$ are not adjacent. Let $\mP_1$ and
$\mP_2$ be two paths in $R$ connecting $u_1$ and $u_2$ such that
$V(\mP_1) \cap V(\mP_2) = \{u_1,u_2\}$ and $|V(\mP_1) \cup V(\mP_2)|$ is
minimal. Hence, the induced subgraph $[V(\mP_i)]$ is a path graph for
$i = 1,2$. Since $R$ is bipartite, both $\mP_1$ and $\mP_2$ have the
same parity.

First assume that $\mP_1$ and $\mP_2$ are even paths. If $v_1 = v_2$,
then we set $\mP_3 := (u_1,v_1,u_2)$ and $\mP_1,\,\mP_2$ and $\mP_3$ are
all even paths connecting $u_1$ and $u_2$; but this contradicts
Theorem \ref{notheta}. If $v_1 \neq v_2$, whenever there exists an
even path $\mP_3'$ in $C$ connecting $v_1, v_2$ we define $\mP_3 :=
(u_1, v_1, \mP_3', v_2, u_2)$, then $\mP_1,\,\mP_2$ and $\mP_3$ are all
even paths connecting $u_1$ and $u_2$; but this is not possible by
Theorem \ref{notheta}. It is easy to check that there exists such an
even path $\mP_3'$ unless if $C$ is an odd partial band consisting of
two odd primitive cycles $C_1 = (a_1,\ldots,a_{r_1},a_1)$ and $C_2 =
(b_1,\ldots,b_{r_2},b_2)$ such that $E(C) = E(C_1) \cup E(C_2) \cup
\{\{a_1,b_1\}\}$ and $\{a_1,b_1\} = \{v_1,v_2\}$. In this situation,
we set $G' := [V(C) \cup V(\mP_1) \cup V(\mP_2)]$ and we have that
$G'$ is a complete intersection and ${\rm deg}_{G'}(a_1) \geq 4$,
then by Corollary \ref{2situaciones} there exists a $v \in V(G')$ of
degree $2$. If $v \in V(C)$, we can assume that $v \in V(C_1)$ and
set $H_1 := [V(C_1) \cup V(\mP_1) \cup \{b_1\}]$, then $b(H_1) =
b(H_1 \setminus \{v\}) = 0$ because $(u_1,\mP_1,u_2,v_2,v_1,u_1)$
 is an odd cycle in $H_1 \setminus \{v\}$. Since $b(C) = b(C \setminus \{v\}) =
 0$, it follows by Lemma \ref{2subgrafos}
that $b(H) = b(H \setminus \{v\})$ where $H =  [V(C) \cap V(H_1)] =
[V(C_1) \cup \{b_1\}]$. Nevertheless, $b(H) = 0$ because $V(C_1)
\subset V(H)$ and $b(H \setminus \{v\}) = 1$ because $H \setminus
\{v\}$ is acyclic, so there is a contradiction. If $v \in V(\mP_1)
\cup V(\mP_2)$, we can assume that $v \in V(\mP_1)$ and we set $H_1 :=
[V(\mP_1)\cup V(C_1) \cup \{b_1\}]$ and $H_2 := [V(\mP_1) \cup
V(\mP_2)]$, then $b(H_1) = b(H_1 \setminus \{v\}) = 0$ and $b(H_2) =
b(H_2 \setminus \{v\}) = 1$. Then by Lemma \ref{2subgrafos}, it
follows that $b(H) = b(H \setminus \{v\})$ where $H =  [V(H_1) \cap
V(H_2)] = [V(\mP_1)]$. Never\-theless, $b(H) = 1$ and $b(H \setminus
\{v\}) = 2$ because $H$ is a path graph and ${\rm deg}_H(v) = 2$, so
there is a contradiction.

Assume now that $\mP_1$ and $\mP_2$ are odd paths. If $v_1 \neq
v_2$, then one can easily find an odd path $\mP_3'$ in $C$
connecting $v_1, v_2$. Therefore, if we set $\mP_3 := (u_1, v_1,
\mP_3', v_2, u_2)$, then $\mP_1, \mP_2$ and $\mP_3$ are all odd
paths connecting $u_1$ and $u_2$; but this is not possible by
Theorem \ref{notheta}. It only suffices to consider the case in
which $v_1 = v_2$. Since $v_1 \in V(C)$, we see that $v_1$ belongs
to the odd primitive cycle $C' = C_1$ or $C' = C_2$ and we set $G'
:= [V(C') \cup V(\mP_1) \cup V(\mP_2)]$. We claim that every vertex
of $C'$ different from $v_1$ has degree $2$. Otherwise there exists
$v' \in V(C'),\ v' \neq v_1$ and $u \in V(\mP_1) \cup V(\mP_2)$ such
that $\{u,v'\} \in E(G')$, we will assume that $u \in V(\mP_1)$.
Then, as we proved before, $\{u,u_1\},\, \{u,u_2\} \in E(R)$. Hence
$(u_1,\mP_2,u_2,u,u_1)$ is an odd cycle in $R$, but this is not
possible because $R$ is bipartite. So we take $v$ any vertex of $C'$
different from $v_1$ and $H_i := [V(\mP_i) \cup V(C')]$ for $i =
1,2$. Then $C_i := (u_1, \mP_i, u_2, v_1, u_1)$ is an odd cycle with
$v \notin V(C_i) \subset V(H_i)$, which gives $b(H_i) = b(H_i
\setminus \{v\}) = 0$. Then by Lemma \ref{2subgrafos}, it follows
that $b(H) = b(H \setminus \{v\})$ where $H =  [V(H_1) \cap V(H_2)]
= [V(C') \cup \{u_1, u_2\}]$. Nevertheless, since $H \setminus
\{v\}$ is acyclic and $V(C) \subset V(H)$, we have that $b(H) = 0$
and $b(H \setminus \{v\}) = 1$, a contradiction.

To sum up, we have proved that whenever $\{u_1,v_1\}, \{u_2,\,v_2\}
\in E(G)$ with $v_1, v_2 \in V(C)$, $u_1, u_2 \in V(R)$ and $u_1
\neq u_2$, then $\{u_1, u_2\} \in E(G)$. If there exist three
different vertices $u_1, u_2, u_3 \in E(G)$ such that $\{u_i,\,v_i\}
\in E(G')$ for some $v_1, v_2, v_3 \in V(C)$, then $u_1, u_2, u_3$
form a triangle in $R$, but this is not possible because $R$ is
bipartite, and the lemma follows. \QED
\end{demo}

\begin{Lemma}\label{12sumas}Let $G$ be a $1$-clique-sum or a $2$-clique-sum of a
graph $H$ and a bipartite ring graph $R$. Then, $G$ is a complete
intersection $\Longleftrightarrow H$ is a complete intersection.
\end{Lemma}
\begin{demo}One implication is obvious because $H$ is an induced subgraph of $G$.
Since bipartite ring graphs are constructed by performing
$1$-clique-sums and $2$-clique-sums of even primitive cycles and
edges, we only have to prove that $G$ is a complete intersection
when it is a $1$-clique-sum or a $2$-clique-sum of a complete
intersection graph $H$ and $K$, where $K$ is either an even
primitive cycle or an edge. If $K$ is an edge $e = \{v_1,v_2\}$ and
$G$ is a $1$-clique-sum of $H$ and $K$, then either ${\rm
deg}_G(v_1) = 1$ or ${\rm deg}_G(v_2) = 1$ and, by Theorem
\ref{principal}, $G$ is a complete intersection. So assume that $K$
is an even primitive cycle $C$. Let $\mathfrak B =
\{B_{w_1},\ldots,B_{w_r}\}$ a minimal set of generators of $P_H$
where $r = {\rm ht}(P_H)$ and consider $\mathfrak B' := \mathfrak B
\cup \{B_C\}$. If we prove that $\mathfrak B'$ generates $P_G$, then
$G$ is a complete intersection because ${\rm ht}(P_G) = {\rm
ht}(P_H) + 1$. We write $B_{w_i} := x^{\alpha_i} - x^{\beta_i}$ for
$i = 1,\ldots,r$ and call $B$ the matrix whose $i$-th row is
$\gamma_i := \alpha_i - \beta_i$, then $B$ is dominating and
$\Delta_r(B) = 1$. We also write $B_{C} := x^{\alpha} - x^{\beta}$
and $B'$ the matrix obtained by adding a new row $\gamma := \alpha -
\beta$ to $B$, let us see that $B'$ is dominating and
$\Delta_{r+1}(B') = 1$. Indeed, $C$ is a cycle which involves at
most one edge of $H$, then by Lemma \ref{dominating} $B'$ is also
dominating and $\gamma$ has only $+1$ and $-1$ in the entries
corresponding to edges in $E(C)$, then $\Delta_{r+1}(B') =
\Delta_r(B) = 1$, which proves the lemma. \QED
\end{demo}

\smallskip
Next we deal with the problems of characterizing when $[V(C) \cup
\{u_1\}]$ is a complete intersection, with $u_1 \in V(R)$, and when
$[V(C) \cup \{u_1,u_2\}]$ is a complete intersection, where $u_1,u_2
\in V(R)$ are adjacent vertices. By Theorem \ref{induce-CI}, when
either $[V(C) \cup \{u_1\}]$ or $[V(C) \cup \{u_1,u_2\}]$ is a
complete intersection, one has that so is $C$ and then, by Theorem
\ref{teoremaestructura}, $C$ is either an odd primitive cycle or an
odd partial band because $C$ is connected. Thus, we will study the
complete intersection property on the following graphs:

\begin{enumerate}
\item $[V(C) \cup \{u_1\}]$, where $C$ is an odd primitive cycle,
\item $[V(C) \cup \{u_1,u_2\}]$, where $u_1,u_2$ are adjacent vertices and $C$ is an odd primitive
cycle,
\item $[V(C) \cup \{u_1\}]$, where $C$ is an odd partial band, and
\item $[V(C) \cup \{u_1,u_2\}]$, where $u_1,u_2$ are adjacent vertices and $C$ is an odd partial band.
\end{enumerate}
The following four lemmas study all these situations. Let us start
with one definition.

\begin{Definition}
An {\it odd partial wheel} $W$ consists of an odd primitive cycle
$C$, a vertex $x \notin V(C)$ and at least one edge connecting $x$
and $C$. The vertex $x$ is called the {\it central vertex of} $W$,
and $C$ is called the {\it principal cycle of} $W$.

Moreover $W$ is a {\it CI-odd-partial-wheel} if $C =
(z_1,\ldots,z_r,z_1)$ and $N_W(x) =\{z_1,z_{s_2}, \ldots,z_{s_k}\}$,
where $k \geq 1$, $1 < s_2 < \cdots < s_{k-1} < s_k$,
$s_3,\ldots,s_k$ are odd and either $s_2 = 2$ or $s_2$ is odd (see
{\rm Figure \ref{opw}}).

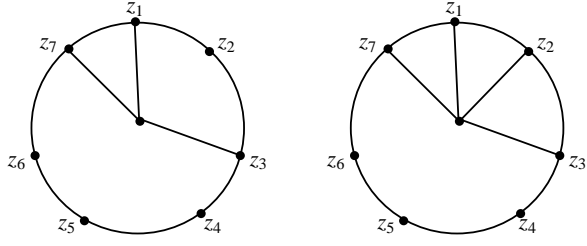
\begin{figure}
\begin{center}
\scalebox{.6} {
\begin{pspicture}(0,-2.385)(12,2.8)
\pscircle[linewidth=0.04,dimen=outer](2.35,-0.035){2.35}
\psdots[dotsize=0.2](2.3,2.315) \psdots[dotsize=0.2](1.18,-2.085)
\psdots[dotsize=0.2](3.74,-1.925) \psdots[dotsize=0.2](4.6,-0.645)
\psdots[dotsize=0.2](3.92,1.655) \psdots[dotsize=0.2](0.84,1.715)
\psdots[dotsize=0.2](0.1,-0.645) \psdots[dotsize=0.2](2.4,0.115)
\psline[linewidth=0.04cm](2.28,2.355)(2.38,0.155)
\psline[linewidth=0.04cm](2.42,0.155)(4.56,-0.625)
\psline[linewidth=0.04cm](2.38,0.135)(0.84,1.675)
\rput(2.3,2.6){\Large $z_1$} \rput(0.8,-2.2){\Large $z_5$}
\rput(4,-2.2){\Large $z_4$} \rput(5,-0.8){\Large $z_3$}
\rput(4.3,1.8){\Large $z_2$} \rput(0.4,1.8){\Large $z_7$}
\rput(-.3,-0.8) {\Large $z_6$}

\pscircle[linewidth=0.04,dimen=outer](9.35,-0.035){2.35}
\psdots[dotsize=0.2](9.3,2.315) \psdots[dotsize=0.2](8.18,-2.085)
\psdots[dotsize=0.2](10.74,-1.925) \psdots[dotsize=0.2](11.6,-0.645)
\psdots[dotsize=0.2](10.92,1.655) \psdots[dotsize=0.2](7.84,1.715)
\psdots[dotsize=0.2](7.1,-0.645) \psdots[dotsize=0.2](9.4,0.115)
\psline[linewidth=0.04cm](9.28,2.355)(9.38,0.155)
\psline[linewidth=0.04cm](9.42,0.155)(11.56,-0.625)
\psline[linewidth=0.04cm](9.38,0.135)(7.84,1.675)
\psline[linewidth=0.04cm](9.4,0.115)(10.9,1.675)
\rput(9.3,2.6){\Large $z_1$} \rput(7.8,-2.2){\Large $z_5$}
\rput(10.9,-2.2){\Large $z_4$} \rput(12,-0.8){\Large $z_3$}
\rput(11.3,1.8){\Large $z_2$} \rput(7.4,1.8){\Large $z_7$}
\rput(6.7,-0.8) {\Large $z_6$}
\end{pspicture}
} \caption{Two CI-odd-partial-wheels}\label{opw}
\end{center}
\end{figure}

\end{Definition}

\begin{Lemma}\label{partialwheel}
Let $W$ be an odd partial wheel. $W$ is a complete intersection if
and only if $W$ is a CI-odd-partial-wheel.
\end{Lemma}
\begin{demo}
Let $x$ be the central vertex and $C$ the principal cycle of $W$, we
denote $r := |V(C)|$.

$(\Rightarrow)$ We proceed by induction on $r$, if $r = 3$ then $W$
is always a CI-odd-partial-wheel. If $r \geq 5$ and ${\rm deg}(x)
\leq 2$ then evidently $W$ is a CI-odd-partial-wheel. So we can
assume that ${\rm deg}_{W}(x) \geq 3$, then there exists $v \in
V(C)$ such that ${\rm deg}_W(v) = 2$, otherwise ${\rm deg}_W(x) = r
> 3$ and this contradicts Corollary \ref{2situaciones}. Thus we
consider $W' := W_v^c$, which is a complete intersection odd partial
wheel and, by induction hypothesis, it is a CI-odd-partial-wheel.
Hence $W'$ has a central vertex $x$ and a principal cycle $C' =
(z_1',\ldots,z_{r-2}')$ such that $N_{W'}(x) =
\{z_1',z_{t_1}',\ldots,z_{t_k}'\}$ where $1 < t_1 < \cdots t_k \leq
r-2$, with $t_2,\ldots,t_k$ odd and $t_1 = 2$ or $t_1$ is odd. Thus
$C = (z_1,\ldots,z_r,z_1)$ and there exists $l \in \{2,\ldots,r-1\}$
such that \begin{center} $\{x,z_i'\} \in E(W') \Longleftrightarrow
\left\{ \begin{array}{cccl} \{x,z_i\} \in E(W) & $ and $ & i \leq
l-1,&$ or $\\ \{x,z_{i+2}\} \in E(W) &$ and $ &i \geq l-1.
\end{array} \right.$
\end{center} If $\{x,z_1\} \notin E(W)$ then we set $y_i := z_{i+2}$ for $1
\leq i \leq r-2$, $y_1 = z_{r-1}$ and $y_2 = z_r$ and have that $W$
is a CI-odd-partial-wheel. If $l \neq 3$ or $\{x,z_4\} \notin E(W)$,
then $W$ also is a CI-odd-partial-wheel. So it remains to study when
$\{x,z_1\},\,\{x,z_4\}$ and $l = 3$, we firstly assume that
$\{x,z_2\} \in E(W)$. If ${\rm deg}_W(x) = 3$, then we set $y_i :=
z_{i+1}$ for $1 \leq i < r$ and $y_r := z_1$ and have that $C =
(y_1,\ldots,y_r,y_1)$ with $N_W(x) = \{y_1,y_3,y_r\}$; thus it is a
CI-odd-partial-wheel. If $N_W(x) = \{z_1,z_2,z_4,z_r\}$, then
setting $y_1 = z_r$, $y_i = z_{i-1}$ for $2 \leq i \leq r$ we have
that $N_W(x) = \{y_1,y_2,y_3,y_5\}$ and $W$ is a
CI-odd-partial-wheel. If ${\rm deg}_W(x) \geq 4$ and $N_W(x) \neq
\{z_1,z_2,z_4,z_r\}$, then we take $j = {\rm min}\{ i >
4\,\vert\,z_i \in N_W(x)\}$ and there are two odd cycles $C_1 :=
(z_1,z_2,x,z_1)$ and $C_2 := (z_4,\ldots,z_j,x,z_4)$ because $j$ is
odd, but there is no edge connecting $C_1$ and $C_2$ and this is
impossible by Lemma \ref{tt}.

Now we assume that $l = 3$ and $\{x,z_2\} \notin E(W)$. If ${\rm
deg}_W(x) = 3$ and $\{x,z_j\} \in E(W)$ for $j = 5$ or $j = r$, then
$W$ is a CI-odd-partial-wheel. Indeed, if $j = r$ we set $y_1 :=
z_r$, $y_i = z_{i-1}$ for $2 \leq i \leq r$, then $N_W(x) =
\{y_1,y_2,y_5\}$ and $W$ is a CI-odd-partial-wheel. If $j = 5$ one
can proceed analogously.

If ${\rm deg}_W(x) \geq 3$ and we are not in the previous
situations, then ${\rm deg}_G(z_2) = 2$ and there exist $5 \leq j_1
\leq j_2 \leq r$ such that $z_{j_1}, z_{j_2} \in N_W(x)$, $(j_1,j_2)
\neq (5,5)$ and $(j_1,j_2) \neq (r,r)$. Since $j_1$ and $j_2$ are
odd, we consider the two even cycles $C_1 :=
(z_1,\ldots,z_{j_1},x,z_1)$ and $C_2 :=
(z_{j_2},\ldots,z_r,z_1,\ldots,z_4,x,z_{j_2})$. Set $H_i :=
[V(C_i)]$ for $i = 1,2$, then $b(H_i) = b(H_i \setminus \{z_2\})$,
however if one takes $H := [V(H_1) \cap V(H_2)]$, then $b(H) = 0$
and $b(H \setminus \{z_2\}) = 1$, which contradicts Lemma
\ref{2subgrafos}.

$(\Leftarrow)$ We write $C = (z_1,\ldots,z_r,z_1)$ and $N_W(x) =
\{z_1,z_{s_2},\ldots,z_{s_k}\}$ where $s_1 := 1 < s_2 < \cdots <
s_k$, $s_3,\ldots,s_k$ are odd and either $s_2 = 2$ or $s_2$ is odd.
If $s_2$ is odd, we set $\overline{W}$ the odd partial wheel with
principal cycle $C' = (z_1',\ldots,z_{r+2}',z_1')$, central vertex
$x'$ and $\{x',z_i'\} \in E(\overline{W})$ if and only if $\{x,z_i\}
\in E(W)$. Clearly ${\rm deg}_{\overline W}(z_{r+1}') = {\rm
deg}_{\overline W}(z_{r+2}') = 2$ and $W =
\overline{W}_{z_{r+1}'}^c$, so if we prove that $\overline{W}$ is a
complete intersection, then by Proposition \ref{contraccion} so is
$W$. We set $R := \overline{W} \setminus \{z_{r+2}'\}$ and have that
$b(\overline{W}) = 0$ and $b(R) = 1$, then by Theorem
\ref{principal}, $\overline{W}$ is a complete intersection if and
only if so is $R$. Since $R$ is a bipartite ring graph, we conclude
that $R$, $\overline{W}$ and $W$ are complete intersections.

Suppose now that $s_2 = 2$, we denote $e_i = \{x,z_{s_i}\}$ for $1
\leq i \leq k$, $e_{k+i} = \{z_i, z_{i+1}\}$ and $e_{k+r} =
\{z_1,z_r\}$. $W$ has $r + 1$ vertices and $r + k$ edges; thus ${\rm
ht}(P_W) = k - 1$. Consider the even cycles $C_1 :=
(x,z_{s_k},z_{s_k+1},\ldots,z_r, z_1,z_2,x)$, $C_2 :=
(x,z_1,z_2,\ldots,z_{s_3},x)$ and $C_i := (x, z_{s_i}, z_{s_i +
1},\ldots,z_{s_{i+1}}, x)$ for $3 \leq i \leq k-1$; then
$$B_{C_1} = x_k  x_{k+1}  x_{k + s_k + 1} \cdots x_{k + r - 1} -
x_2  x_{k + s_k} \cdots x_{k + r},$$ $$B_{C_2} = x_1 x_{k+2} \cdots
x_{k + s_3 - 1}  - x_3 x_{k+1} \cdots x_{k + s_3 - 2}, {\rm \ and}$$
$$B_{C_i} = x_i  x_{k + s_i + 1} \cdots x_{k + s_{i+1} - 1} -
x_{i+1}  x_{k + s_i} \cdots x_{k + s_{i+1} - 2} {\rm \ for \ } 3
\leq i \leq k-1,$$  let us prove that $P_W =
(B_{C_1},\ldots,B_{C_{k-1}}).$ We set $B_{C_i} = x^{\alpha_i} -
x^{\beta_i}$ and $\gamma_i := \alpha_i - \beta_i$ for $1 \leq i \leq
k-1$; then, $$\gamma_1 := - e_2 + e_k + e_{k+1} - e_{k + s_k} + e_{k
+ s_k + 1} - \cdots + e_{k + r - 1} - e_{k + r} \in \Z^{k+r},$$
$$\gamma_2 := e_1 - e_3 - e_{k+1} + e_{k+2} - \cdots - e_{k + s_3 -
2} + e_{k + s_3 - 1} \in \Z^{k+r}, {\rm \ and}$$
$$\gamma_i := e_i - e_{i+1} - e_{k + s_i} + e_{k+s_i+1} - e_{k+s_i+2} + \cdots -
e_{k + s_{i+1} - 2} + e_{k + s_{i+1} -1} \in \Z^{k + r}$$  for $3
\leq i \leq k-1$, and denote by $B$ the $(k - 1) \times (k + r)$
matrix whose $i$-th row is $\gamma_i$. It is evident that
$\Delta_{k-1}(B) = 1$ and for every $j \in \{1,2,k+2,\ldots,k+r\}$
the $j$-th column of $B$ has only one nonzero entry; thus by Lemma
\ref{dominating} $B$ is dominating if and only if $B'$ is dominating
where $B'$ is the $k-1 \times k-1$ matrix consisting of the columns
$3,4,\ldots,k+1$ of $B$. $B'$ has exactly two nonzero entries in
each row and in each column and both nonzero entries in the first
row of $B'$ are positive. Hence, $B'$ is dominating. Therefore we
conclude that $W$ is a complete intersection and $P_W =
(B_{C_1},\ldots,B_{C_{k-1}})$.
 \QED
\end{demo}

\begin{Definition}A connected graph $G$ is called a {\it CI-double-wheel}
if  its vertex set is $V(G) = V(C) \cup \{b_1,b_2\}$, where $C =
(a_1,\ldots,a_r,a_1)$ is an odd primitive cycle and $E(G) = E(C)
\cup \left\{ \{b_1,b_2\}, \{b_1,a_{j_1}\}, \ldots,\{b_1,a_{j_s}\},
\{b_2,a_{k_1}\},\ldots, \{b_2,a_{k_t}\} \right\},$ for some $s,t
\geq 1$, $1 \leq j_1 < \cdots < j_s \leq k_1 < \cdots < k_t \leq r$
and $j_1,\ldots,j_s,k_1,\ldots,k_t$ are odd (see {\rm Figure
\ref{doubwheel}}).

\begin{figure}
\begin{center}
\scalebox{.8} {
\begin{pspicture}(0,-2)(3.7,1.8)
\psarc[linewidth=0.03](1.8,0.085){1.73}{0.0}{180.0}
\psline[linewidth=0.03cm](0.09,0.095)(3.51,0.135)
\psdots[dotsize=0.14](0.07,0.075) \psdots[dotsize=0.14](0.71,0.115)
\psdots[dotsize=0.14](1.31,0.135) \psdots[dotsize=0.14](1.87,0.135)
\psdots[dotsize=0.14](2.51,0.155) \psdots[dotsize=0.14](2.99,0.155)
\psdots[dotsize=0.14](3.51,0.155) \psdots[dotsize=0.14](0.89,-1.745)
\psdots[dotsize=0.14](2.53,-1.705)
\psline[linewidth=0.03cm](0.89,-1.725)(2.49,-1.685)
\psline[linewidth=0.03cm](0.05,0.075)(0.87,-1.765)
\psline[linewidth=0.03cm](1.27,0.155)(0.91,-1.725)
\psline[linewidth=0.03cm](1.33,0.135)(2.51,-1.705)
\psline[linewidth=0.03cm](2.53,-1.685)(3.47,0.215)
\rput(-0.2,.3){$a_1$} \rput(0.7,.4){$a_2$} \rput(1.3,.4){$a_3$}
\rput(1.85,.4){$a_4$} \rput(2.5,.4){$a_5$} \rput(3,.4){$a_6$}
\rput(3.7,.3){$a_7$} \rput(0.9,-2){$b_1$} \rput(2.5,-2){$b_2$}
\end{pspicture}
} \caption{A CI-double-wheel}\label{doubwheel}
\end{center}
\end{figure}
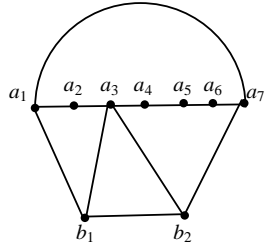

\end{Definition}

\begin{Lemma}\label{doublewheel}Let $G$ be a connected graph with $V(G) = V(C) \cup \{b_1,b_2\}$
where $C$ is an odd primitive cycle, $\{b_1, b_2\} \in E(G)$ and
${\rm deg}_G(b_1), {\rm deg}_G(b_2) \geq 2$. Then, $G$ is a complete
intersection if and only if $G$ is a CI-double-wheel.
\end{Lemma}
\begin{demo}($\Rightarrow)$ We proceed by induction on $r := |V(C)|$. If $r = 3$, then by
Corollary \ref{cotasuperior} it has $\leq 7$ edges. If every vertex
of $C$ has degree $\geq 3$, then $G$ has a subgraph $\mK_{2,3}$,
which contradicts Corollary \ref{subgraphK23}. Thus one can write
$V(G) = V(C) \cup \{b_1,b_2\}$, where $C = (a_1,a_2,a_3,a_1)$, ${\rm
deg}_G(a_2) = {\rm deg}_G(b_2) = 2$ and ${\rm deg}_G(b_1) \leq 3$,
thus $G$ is a CI-double wheel. Assume now that $r \geq 5$, we claim
that there exists a vertex of degree $2$ in $C$. Indeed, if ${\rm
deg}_G(v) \geq 3$ for every $v \in V(C)$, then we have that $4r + 2
\leq 2 |E(G)| \leq 3 |V(G)| = 3r + 6$, which contradicts Corollary
\ref{cotasuperior}. Therefore we take $v \in V(C)$ of degree $2$ and
consider $G' := G_v^c$, which is a CI-double wheel by induction
hypothesis. Thus $V(G') = V(C') \cup \{b_1,b_2\}$, where $C' =
(a_1',\ldots,a_{r-2}',a_1')$ is an odd primitive cycle and $E(G') =
E(C') \cup \left\{ \{b_1,b_2\},  \{b_1, a_{j_1}'\}, \ldots, \{b_1,
a_{j_s}'\}, \{b_2,a_{k_1}'\},\ldots,\{b_2,a_{k_t}'\}\right\}$ with
$1 \leq j_1 < \cdots < j_s \leq k_1 < \cdots < k_t \leq r$ and
$j_1,\ldots,j_s,k_1,\ldots,k_t$ are odd. Moreover, there exists $l
\in \{2,\ldots,r-1\}$ such that
\begin{center} $\{b_i,a_j'\}
\in E(G') \Longleftrightarrow \left\{ \begin{array}{cccl}
\{b_i,a_j\} \in E(G) & $ and $ & j \leq l - 1,&$ or $\\
\{b_i,a_{j+2}\} \in E(G) &$ and $ &j \geq l - 1.
\end{array} \right.$
\end{center}
If $\{b_1,a_{l+1}\} \notin E(G)$, $\{b_2,a_{l-1}\} \notin E(G)$ or
${\rm deg}(b_1) = {\rm deg}(b_2) = 2$, then $G$ is a CI-double
wheel. Therefore if $G$ is not a CI-double-wheel one can assume that
$\{b_1,a_{l+1}\},  \{b_2,a_{l-1}\}  \in E(G)$ and ${\rm deg}(b_1) >
2$. Thus, $l - 1 = j_s = k_1,\ldots,k_t$ are odd and
$\{b_1,a_{j_1}\},  \{b_1,a_{l + 1}\},  \{b_2,a_{l - 1}\} \in E(G)$.

We separate two cases, if $b_2$ has degree $> 2$ then
$\{b_2,a_{k_t}'\} \in E(G')$ and $\{b_2,a_{k_t + 2}\} \in E(G)$,
then there exist three even paths
$$\mP_1 := (a_{l-1},a_l,a_{l+1}),\,  \mP_2 :=
(a_{l-1},a_{l-2},\ldots,a_{j_1},b_1,a_{l+1}) {\rm\ and}$$ $$\mP_3 :=
(a_{l-1},b_2,a_{k_t+2},a_{k_t+1},\ldots,a_{l+1})$$ connecting
$a_{l-1}$ and $a_{l+1}$ and $V(\mP_i) \cap V(\mP_j) =
\{a_{l-1},a_{l+1}\}$ for $1 \leq i < j \leq 3$, which is not
possible by Theorem \ref{notheta}. If $b_2$ has degree $2$ and $j_1
= l-1$, then $G$ is a CI-double-wheel. Hence we assume that $b_2$
has degree $2$ and $j_1 < l-1$ and we consider the even cycles $C_1
:= (b_1,b_2,a_{l-1},a_l,\ldots,a_r,a_1,\ldots,a_{j_1},b_1)$ and $C_2
:= (b_1,b_2,a_{l-1},a_{l-2},\ldots,a_1,a_r,\ldots,a_{l+1},b_1)$, and
set $H_i := [V(C_i)]$, then $b(H_i) = b(H_i \setminus \{b_2\}) = 1$.
However setting $H := [V(C_1) \cap V(C_2)] =
[\{a_1,\ldots,a_{j_1},a_{l-1},a_{l+1},\ldots,a_r,b_1,b_2\}]$, if
$\{b_1,a_{l-1}\} \notin E(G)$, then $b(H) = 0$ because $H$ is
connected and the odd cycle $C_3 :=
(b_1,a_{l+1},\ldots,a_r,a_1,\ldots,a_{j_1},b_1)$ is in $H$, and $b(H
\setminus \{b_2\}) = 1$ because $a_{l-1}$ is an isolated vertex in
$H \setminus \{b_2\}$. It only remains to consider the situation in
which $\{b_1,a_{l-1}\} \in E(G)$; in this case $G$ is a
$2$-connected graph, there are two odd cycles $C_3 =
(b_1,a_{l-1},b_2,b_1)$ and $C_4 =
(b_1,a_{l+1},\ldots,a_r,a_1,\ldots,a_{j_1},b_1)$ with a vertex in
common and there is no edge connecting them, a contradiction to
Lemma \ref{tt}.

$(\Leftarrow)$ Denote by $G'$ the graph obtained by adding a new
vertex $b_3$ and two edges $\{b_1,b_3\}$ and $\{b_2,b_3\}$, then
$G'$ is an odd partial band and by Proposition \ref{oddpartilband}
$G'$ is a complete intersection. Furthermore, $G = G' \setminus
\{b_3\}$, then $G$ is a complete intersection. \QED
\end{demo}

\begin{Definition}A {\it CI-vertex-band} is a graph $G$
with vertices $V(C_1) \cup V(C_2) \cup \{c\}$, where $C_1 =
(a_1,\ldots,a_r,a_1)$ and $C_2 = (b_1,\ldots,b_s,b_1)$ are vertex
disjoint odd primitive cycles and $$E(G) = E(C_1) \cup E(C_2) \cup
\left\{ \{a_1,b_1\},\{a_1,b_{i_2}\} \ldots \{a_1,b_{i_k}\}, \{c,
a_2\}, \{c, a_r\} \right\},$$ for some $k \geq 1$, $i_2 < \cdots <
i_k \leq s$ and $i_2,\ldots,i_k$ are odd (see {\rm Figure
\ref{vertexb}}).

\begin{figure}
\begin{center}
\scalebox{.8} {
\begin{pspicture}(0.5,-1.8)(7.3,1.7)
\pscircle[linewidth=0.04,dimen=outer](1.7371875,-0.1328125){1.53}
\pscircle[linewidth=0.04,dimen=outer](5.6371875,-0.1528125){1.45}
\psdots[dotsize=0.14](3.2671876,-0.0828125)
\psdots[dotsize=0.14](2.0671875,1.3371875)
\psdots[dotsize=0.14](0.3671875,0.5571875)
\psdots[dotsize=0.14](0.6671875,-1.2228125)
\psdots[dotsize=0.14](2.1071875,-1.5628124)
\psdots[dotsize=0.14](4.1871877,-0.0828125)
\psdots[dotsize=0.14](4.4871874,0.6371875)
\psdots[dotsize=0.14](5.0271873,1.1371875)
\psdots[dotsize=0.14](6.1471877,1.1571875)
\psdots[dotsize=0.14](7.0471873,-0.1828125)
\psdots[dotsize=0.14](6.3871875,-1.3428125)
\psdots[dotsize=0.14](4.7071877,-1.2228125)
\psdots[dotsize=0.14](1.4671875,-0.1028125)
\psline[linewidth=0.04cm](2.0271876,1.3571875)(1.4671875,-0.0428125)
\psline[linewidth=0.04cm](1.4671875,-0.1428125)(2.1071875,-1.5828125)
\psline[linewidth=0.04cm](3.2871876,-0.0628125)(4.1671877,-0.0828125)
\psline[linewidth=0.04cm](3.2671876,-0.0828125)(4.9671874,1.1771874)
\psline[linewidth=0.04cm](3.2471876,-0.0828125)(4.6871877,-1.2028126)
\rput(3.0,0){$a_1$} \rput(2.1732812,1.5){$a_2$}
\rput(0.22578125,.8){$a_3$} \rput(0.41453126,-1.4){$a_4$}
\rput(2.2889063,-1.8){$a_5$} \rput(4.5,-0.0328125){$b_1$}
\rput(4.8,0.7071875){$b_2$} \rput(5.1626563,1.5){$b_3$}
\rput(6.6314063,1.2){$b_4$} \rput(7.3,-0.0328125){$b_5$}
\rput(6.730781,-1.4328125){$b_6$} \rput(4.6296873,-1.4128125){$b_7$}
\rput(1.3751563,0.1){$c$}
\end{pspicture}
} \caption{A CI-vertex-band}\label{vertexb}
\end{center}
\end{figure}
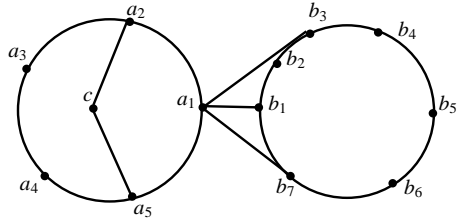
\end{Definition}

\begin{Lemma}\label{civertexband}Let $G$ be a connected graph with $V(G) = V(C) \cup
\{c\}$ where $C$ is an odd partial band. Then, $G$ is a complete
intersection if and only if ${\rm deg}_G(c) = 1$ or $G$ is a
CI-vertex-band.
\end{Lemma}
\begin{demo}($\Rightarrow$) We denote by $C_1$ and $C_2$ the two vertex
disjoint odd primitive cycles such that $V(C) = V(C_1) \cup V(C_2)$.
We first aim to prove that $c$ has degree $\leq 2$. For this purpose
we will prove that if $v \in V(C)$ has degree $2$ and does not
belong to a triangle, then ${\rm deg}_G(c) = {\rm deg}_{G_v^c}(c)$.
Suppose that ${\rm deg}_G(c) > {\rm deg}_{G_v^c}(c)$, this means
that $N_G(v) = \{v_1,v_2\}$ and $\{v_1,c\}, \{v_2,c\} \in E(G)$.
Then we consider $H := [\{v,v_1,v_2,c\}]$ and have that $v \in V(H)
\cap V(C)$, $b(H) = b(H \setminus \{v\}) = 1$, $b(C) = b(C \setminus
\{v\}) = 0$. However, considering $H' := [V(H) \cap V(C)] =
[\{v,v_1,v_2\}]$, then $b(H' \setminus \{v\}) = b(H') + 1$, which
contradicts Lemma \ref{2subgrafos}.

Suppose that ${\rm deg}_G(c) \geq 3$. For every $v \in V(C)$ of
degree $2$ which does not belong to a triangle we consider $G_c^v$
and we repeat this until we get a graph $G'$ in which every vertex
of degree $2$ belongs to a triangle, then we have proved that ${\rm
deg}_G(c) = {\rm deg}_{G'}(c)$. Note that $V(G') = V(C') \cup \{c\}$
where $C'$ is an odd partial wheel with primitive cycles $C_1'$ and
$C_2'$. Since $G'$ has an odd number of vertices, we have that it
cannot be a band or a M\"obius band and there exists a $v \in V(G')$
of degree $2$, say $v \in V(C_1')$. By construction $C_1'$ has to be
a triangle, let us prove that $C_2'$ is also a triangle. Suppose
that $C_2'$ is not a triangle, then $s := |V(C_2')| \geq 5$ and
${\rm deg}_{G'}(u) \geq 3$ for every $u \in V(C_2')$, which implies
that $|E(G')| \geq 2s + 3$ and $|V(G')| = s + 4$ and by Corollary
\ref{cotasuperior} we have that $2 |E(G')| \leq 3 |V(G')|$; thus $s
= 5$ and $|E(G')| = 13$. This means that ${\rm deg}_{G'}(z) = 3$ for
every $z \in V(C_2')$ and $\{c,u\} \notin E(G')$ for every $u \in
V(C_1')$. By symmetry, one can assume that $C_1' = (a_1,a_2,a_3)$,
$C_2' = (b_1,b_2,b_3,b_4,b_5,b_1)$ and $E(G') = E(C_1') \cup E(C_2')
\cup \{\{a_1,b_1\}\} \cup E$, where $E$ is one of these:
\begin{itemize}[leftmargin=.8cm]
\item[{\rm (a)}] $E = \{\{b_2,c\}, \{b_3,c\}, \{b_4,c\}, \{b_5,c\} \},$
\item[{\rm (b)}] $E = \{\{b_2,a_1\}, \{b_3,c\}, \{b_4,c\}, \{b_5,c\} \},$
\item[{\rm (c)}] $E = \{\{b_2,a_2\}, \{b_3,c\}, \{b_4,c\}, \{b_5,c\} \},$
\item[{\rm (d)}] $E = \{\{b_2,c\}, \{b_3,a_1\}, \{b_4,c\}, \{b_5,c\} \},$ or
\item[{\rm (e)}] $E = \{\{b_2,c\}, \{b_3,a_2\}, \{b_4,c\}, \{b_5,c\} \}.$
\end{itemize}
If (a) occurs we set $H_i := [V(C_1') \cup \{c,b_1,b_{i-1},b_i\}]$
and have that $b(H_i) = b(H_i \setminus \{a_2\}) = 0$ for $i = 3,5$;
however if one sets $H := [V(H_3) \cap V(H_5)] = [V(C_1') \cup
\{c,b_1\}]$, then $b(H) = 1$ and $b(H \setminus \{a_2\}) = 2$, which
is impossible by Lemma \ref{2subgrafos}. If (b) or (c) holds, we set
$\mathcal P_1 := (b_3,b_4,b_5)$, $\mathcal P_2 := (b_3,c,b_5)$, in
(b) we also set $\mathcal P_3 := (b_3,b_2,a_1,b_1,b_5)$ and in (c)
we also set $\mathcal P_3 := (b_3,b_2,a_2,a_3,a_1,b_1,b_5)$; in both
situations we have three even paths connecting $b_3$ and $b_5$, but
this is not possible by Theorem \ref{notheta}. In the last two cases
we set $\mathcal P_1 := (b_3,b_2,b_1)$, $\mathcal P_2 :=
(b_3,b_4,c,b_5,b_1)$,  in (d) we also set $\mathcal P_3 :=
(b_3,a_1,b_1)$ and in (e) we also set $\mathcal P_3 :=
(b_3,a_2,a_3,a_1,b_1)$; but this is not possible again by Theorem
\ref{notheta}.

So assume now that both $C_1'$ and $C_2'$ are triangles and ${\rm
deg}_{G'}(c) \geq 3$, then necessarily ${\rm deg}_{G'}(c) = 3$ and
$E(C') = E(C_1') \cup E(C_2') \cup \{\{a_1,b_1\}\}$, otherwise $2
|E(G')| > 3 |V(G')|$. Then one can write $C_1' = (a_1,a_2,a_3,a_1)$,
$C_2' = (b_1,b_2,b_3,b_1)$ and $N_{G'}(c)$ is one of these:
\begin{itemize}
\item $N_{G'}(c) = \{b_1,b_2,b_3\},$
\item $N_{G'}(c) = \{a_1,b_1,b_2\},$
\item $N_{G'}(c) = \{a_2,b_1,b_2\},$
\item $N_{G'}(c) = \{a_1,b_2,b_3\}$, or
\item $N_{G'}(c) = \{a_2,b_2,b_3\}.$
\end{itemize}

We set $u := a_3$, $H_1 := C'$ and $H_2 := G' \setminus \{b_3\}$ in
the first three cases and $H_2 := G' \setminus \{b_1\}$ in the last
two. In all of them $b(H_1) = b(H_1 \setminus \{u\}) = 0$, $b(H_2) =
b(H_2 \setminus \{u\}) = 0$. However, $b(H) \neq b(H \setminus
\{u\})$ where $H := [V(H_1) \cap V(H_2)]$, a contradiction to Lemma
\ref{2subgrafos}. Thus ${\rm deg}_G(c) \leq 2$.

If ${\rm deg}_G(c) = 2$, since $b(G) = b(G \setminus \{c\}) = 0$, by
Theorem \ref{principal} we get that $P_G = P_{C} \cdot
k[x_1,\ldots,x_n] + (B_w)$, where $w$ is an even closed walk with
$$V(w) = \{c\} \cup N_G(c) \cup \{v \in V(G) \, \vert \, b(G
\setminus \{v\}) < b(C \setminus \{v\})\}.$$ We assume that $C_1 =
(a_1,\ldots,a_r,a_1)$ and $C_2 = (b_1,\ldots,b_s,b_1)$ and that
$\{a_1,b_1\} \in E(C)$. Moreover, if $v \in V(w)$ and $v \notin
N_G(c) \cup \{c\}$, then $v \in \{a_1,b_1\}$ because $b(C \setminus
\{a_j\}) = b(C \setminus \{b_j\}) = 0$ for all $j \geq 2$. Therefore
we can suppose that $a_1 \in V(w)$ and $a_1 \notin N_G(c)$; thus
$b(G \setminus \{a_1\}) = 0$ and $b(C \setminus \{a_1\}) = 1$. Since
$b(C \setminus \{a_1\}) = 1$, we can assume that $E(C) = E(C_1) \cup
E(C_2) \cup \{\{a_1,b_1\},\{a_1,b_{j_2}\},\ldots,\{a_1,b_{j_k}\}\}$
for some $k \geq 1$, and $j_2,\ldots,j_k$ are odd because $C$ is an
odd partial band. Moreover, since $b(G \setminus \{a_1\}) = 0$ it
follows that $N_G(c) = \{a_i,b_j\}$ for some $1 < i \leq r$, $1 \leq
j \leq s$ or $N_G(c) = \{a_i,a_j\}$ for some $1 < i < j \leq r$ and
$i \not\equiv j \ ({\rm mod}\ 2)$. If $N_G(c) = \{a_i,a_j\}$, then
$b_1 \notin V(w)$ because $b(G \setminus \{b_1\}) = b(C \setminus
\{b_1\})$. Additionally, if $N_G(c) = \{a_i,b_j\}$ for some $1 < i
\leq r$, $1 \leq j \leq s$, since $G$ can not be $2$-connected by
Lemma \ref{tt}, we get that $j = 1$ and $\{a_1,b_1\}$ is the only
edge connecting $C_1$ and $C_2$. Putting all together, we can assume
that one of these occurs:
\begin{itemize}[leftmargin=1cm]
\item[{\rm (a)}] $V(w) = \{c,a_1,a_i,b_1\}$, where
$N_G(c) = \{a_i,b_1\}$ with $1 < i \leq r$ and $E(C) = E(C_1) \cup
E(C_2) \cup \{\{a_1,b_1\}\}$, or
\item[{\rm (b)}] $V(w) = \{c,a_1,a_i,a_j\}$, where
$N_G(c) = \{a_i,a_j\}$ with $1 < i < j \leq r$,\, $i \not\equiv j\
({\rm mod}\ 2)$ and $E(C) = E(C_1) \cup E(C_2) \cup
\{\{a_1,b_1\},\{a_1,b_{j_2}\},\ldots,\{a_1,b_{j_k}\}\}$ for some $k
\geq 1$, and $j_2,\ldots,j_k$ are odd.
\end{itemize}

In both cases $w$ is a length $4$ cycle by Lemma \ref{primitivo}. In
(a) we have that $w = (c, a_i, a_1, b_1)$, and we can assume that $i
= 2$. We proved in Proposition \ref{oddpartilband} that $P_C =
(B_{w'})$ where $w' = (a_1,\ldots,a_r,a_1,b_1,\ldots,b_s,b_1,a_1)$;
thus $P_G = (B_{w},B_{w'})$. However this is not possible because
denoting $e_1$ and $e_2$ the edges $\{a_1,b_1\}$ and $\{a_1,a_2\}$,
then $(B_{w},B_{w'}) \subsetneq J := (x_1,x_2)$ and $2 = {\rm
ht}(P_G) < {\rm ht}(J) = 2$. Finally, if (b) occurs we have that $w
= (c, a_i, a_1, a_j)$, which  implies that $i = 2$, $j = r$.
Therefore $G$ is a CI-vertex-band.

$(\Leftarrow)$ If ${\rm deg}_G(c) = 1$, $G$ is a complete
intersection if and only if so is $C$ and $C$ is an odd partial
band, which is a complete intersection by Proposition
\ref{oddpartilband}.

If $G$ is a CI-vertex-band with $V(G) = V(C) \cup \{c\}$, where $C$
 consists of two odd vertex disjoint cycles $C_1
= (a_1,\ldots,a_r,a_1)$ and $C_2 = (b_1,\ldots,b_s,b_1)$ and $E(G) =
E(C_1) \cup E(C_2) \cup \left\{ \{a_1,b_{i_1}\}, \ldots
,\{a_1,b_{i_k}\}, \{c, a_2\}, \{c, a_r\} \right\}$ where $1 =
b_{i_1} < \cdots < b_{i_k}$ and $b_{i_1},\ldots,b_{i_k}$ are odd.
Let us prove that $P_G = P_C \cdot k[x_1,\ldots,x_n] + (B_w)$ where
$w = (c,a_r,a_1,a_2,c)$.

We have that ${\rm ht}(P_G) = k + 1$, we set $e_j := \{a_{1},
b_{i_j}\}$  for $1 \leq i \leq k$, $e_{k+j} := \{a_j, a_{j+1}\}$ for
$1 \leq j < r$, $e_{k + r} := \{a_1, a_{r}\}$, $e_{k + r + j} :=
\{b_j, b_{j + 1}\}$ for $1 \leq j < s$ and $e_{k + r + s} := \{b_1,
b_{s}\}$, $e_{k+r+s+1} := \{c, a_2\}$ and $e_{k+r+s+2} := \{c,
a_r\}$.

For every $j \in \{1,\ldots,k-1\}$, let $w_j$ be the even primitive
cycle
$$w_j := (b_{i_j}, a_{1}, b_{i_{j+1}},  b_{i_{j+1} - 1},\ldots,
b_{i_j})$$ and $w_k := (b_{i_k}, a_1,C_1, a_1, b_{1}, b_{s},\ldots,
b_{i_k})$.

If we denote $B_{w_j} = x^{\alpha_j} - x^{\beta_j}$ with $\alpha_j,
\beta_j \in \N^{k + r + s + 2}$ for every $j \in \{1,\ldots,k\}$ and
we let $B$ be the $k \times (k + r + s + 2)$ matrix whose $j$-th row
is $\alpha_j - \beta_j$; then we proved in Proposition
\ref{oddpartilband} that $B$ is dominating and $\Delta_k(B) = 1$.

We denote by $B'$ the matrix obtained by adding to $B$ the row
$e_{k+1} - e_{k + r} - e_{k+r+s+1} + e_{k+r+s+2} \in \Z^{k+r+s+2}$.
Then $\Delta_{k+1}(B') = 1$ because $\Delta_k(B) = 1$, if we prove
that $B'$ is dominating then $P_G = (B_{w_1},\ldots,B_{w_k}, B_{w})$
and it is a complete intersection. Since the columns $k+r+s+1$ and
$k+r+s+2$ have only one nonzero entry, if we denote by $B''$ the
matrix obtained by removing these two columns from $B'$, by Lemma
\ref{dominating} we get that $B'$ is dominating if and only if so is
$B''$. Assume that $B''$ has a mixed square submatrix $D$, since $B$
is dominating then the last row of $B''$ is in $D$, moreover
 the columns $k+1$ and $k+r$ of $B''$ have to be in $D$ because
these are the only two nonzero entries in the last row of $B''$.
Furthermore the columns $k+1$ and $k+r$ of $B''$ have only two
nonzero entries, which are those in the rows $k$ and $k+1$ and the
entries in the row $k$ are both negative. So, if we remove the last
row of $D$ and the column $k+1$ we get $D'$ another square matrix of
$B''$ which is mixed; but $D'$ is also a submatrix of $B$; which it
is a contradiction. \QED
\end{demo}

\begin{Lemma}\label{2oddring}Let $G$ be a connected graph with $V(G) = V(C) \cup
\{c_1,c_2\}$ where $C$ is an odd partial band, $\{c_1,c_2\} \in
E(G)$ and $c_1, c_2$ have degree $\geq 2$. Then, $G$ is a complete
intersection if and only if $G$ is a $2$-clique-sum of $C$ and a
length $4$ cycle $C$.
\end{Lemma}
\begin{demo}$(\Rightarrow)$ For $i = 1,2$ we have that $G_i := G \setminus \{c_i\}$ is a
complete intersection, then by Lemma \ref{civertexband} either ${\rm
deg}_{G_i}(c_{3-i}) = 1$ or $G_i$ is a CI-vertex-band. Firstly
assume that $G_1$ and $G_2$ are CI-vertex-bands and denote $C_1 =
(a_1,\ldots,a_r,a_1)$ and $C_2 = (b_1,\ldots,b_s,b_1)$ with
$\{a_1,b_1\} \in E(G)$ the two odd primitive cycles such that $V(C)
= V(C_1) \cup V(C_2)$.
 If $\{c_i,a_j\} \in E(G)$ for $i = 1,2$, $j
= 2,r$, then $G$ has a subgraph $\mathcal K_{2,3}$ with vertices
$\{c_1,c_2,a_1,a_2,a_r\}$ and this is not possible. If one has that
$\{c_1,a_r\}, \{c_1,a_2\}, \{c_2,b_s\}, \{c_2,b_2\} \in E(G)$, then
there are three even paths connecting $c_1$ and $a_1$, namely
$\mathcal P_1 := (c_1,a_2,a_1)$, $\mathcal P_2 := (c_1,a_r,a_1)$ and
$\mathcal P_3 := (c_1,c_2,b_2,b_1,a_1)$ and this is not possible by
Theorem \ref{notheta}. So assume now that ${\rm deg}_G(c_1) = 2$ and
$G_1$ is a CI-vertex-band; then $P_G = P_{G_1} \cdot
k[x_1,\ldots,x_n] + (B_w)$ where $w$ is an even closed walk with
$V(w) = \{c_1\} \cup N_G(c_1) \cup \{u \in V(G)\, \vert \, b(G
\setminus \{u\}) < b(G \setminus \{u,c_1\})\}$. Since $c_1, c_2 \in
V(w)$, we get that $a_2$ or $a_r \in V(w)$, say $a_2 \in V(w)$. But
$b(G \setminus \{c_1, a_2\}) = 0$, so $a_2 \in N_G(c_1)$ and $V(w) =
\{c_1,c_2,a_2,b_1\}$, but such a closed walk does not exist.

Then we have proved that ${\rm deg}_G(c_1) = {\rm deg}_G(c_2) = 2$
and there exist $u_1, u_2 \in V(C)$ such that $\{c_1,u_1\},
\{c_2,u_2\} \in E(G)$. Then by Theorem \ref{principal}, $P_G = P_{G
\setminus \{c_2\}} \cdot k[x_1,\ldots,x_n] + (B_w)$ where $w$ is an
even closed walk with $V(w) = \{c_2\} \cup N_G(c_2) \cup \{u \in
V(G)\, \vert \, b(G \setminus \{u\}) < b(G \setminus \{u,c_2\})\}$.
Then, $w = (u_1,c_1,c_2,u_2,v_1,\ldots,v_t,u_1)$ for some $t \geq 0$
and $b(G \setminus \{v_i\}) < b(G \setminus \{v_i,c_2\})$ for $1
\leq i \leq t$. Nevertheless, $b(G \setminus \{z,c_2\}) \geq 1$ if
and only if $z = u_1$ or $z = a_1$ and ${\rm deg}_C(a_i) = 2$ for
every $i \geq 2$ or $z = b_1$ and ${\rm deg}_C(b_i) = 2$ for every
$i \geq 2$. Then $\{c_1,c_2,u_1,u_2\} \subset V(w) \subset
\{c_1,c_2,u_1,u_2,a_1,b_1\}$. If $|V(w)| = 6$, then we can assume
that $u_1 = a_2$, $u_2 = b_2$ and $\{a_1,b_1\}$ is the only edge
connecting $C_1$ and $C_2$ but this contradicts Lemma \ref{tt}. If
$|V(w)| = 5$, then $w$ is not an even cycle and by Lemma
\ref{primitivo}, $u_1 = u_2$; however, if $u_1 \in V(C_1)$, then
$b(G \setminus \{b_1,c_2\}) = b(G \setminus \{b_1\})$ and $b_1
\notin V(w)$, which is a contradiction. Then $|V(w)| = 4$, $w$ is a
cycle and $\{u_1,u_2\} \in E(G)$; so $G$ is a $2$-clique-sum of $C$
and the length $4$ cycle $(c_1,c_2,u_2,u_1,c_1)$.

$(\Leftarrow)$ It follows directly from Lemma \ref{12sumas} and
Proposition \ref{oddpartilband}.
 \QED
\end{demo}

\bigskip

Now we can state and prove the following characterization theorem,
which allows us to list all families of complete intersection graphs
$G = [C;R]$ such that $R$ is $2$-connected and $C$ is connected.

\begin{Theorem}\label{ultimo}Let $G = [C; R]$ be a connected graph. If $R$ is $2$-connected
and $C$ is connected, then $G$ is a complete intersection if and
only if $G$ is one of the following graphs:
\begin{itemize}
\item[{\rm (a)}] a bipartite ring graph,
\item[{\rm (b)}] a $1$-clique-sum of a bipartite ring graph and either
\begin{itemize} \item[{\rm (b.1)}] a CI-odd-partial-wheel, \item[{\rm (b.2)}] a $1$-clique-sum of an odd partial band and an
edge, or
\item[{\rm (b.3)}] a CI-vertex-band,
\end{itemize}
\item[{\rm (c)}] a $2$-clique-sum of a bipartite ring graph and either \begin{itemize}  \item[{\rm (c.1)}] a CI-double-wheel, or
\item[{\rm (c.2)}]an odd partial band. \end{itemize}
\end{itemize}
\end{Theorem}
\begin{demo}$(\Rightarrow)$ Since $G$ is connected, by Theorem \ref{teoremaestructura}
we have that $R$ is a bipartite ring graph and $C$ is either the
empty graph, an odd primitive cycle or an odd partial band. If $C$
is the empty graph, then $G = R$ is a bipartite ring graph.
Otherwise, by Proposition \ref{propestructura}, either
\begin{itemize} \item[(1)] there exists a $b_1 \in V(R)$ such that $[V(C) \cup
\{b_1\}]$ is a complete intersection, or \item[(2)] there exist two
adjacent vertices $b_1,b_2 \in V(R)$ such that $[V(C) \cup
\{b_1,b_2\}]$ is a complete intersection. \end{itemize} Assume (1)
holds, if $C$ is an odd primitive cycle, by Lemma \ref{partialwheel}
we obtain (b.1), and if $C$ is an odd partial band, by Lemma
\ref{civertexband} we obtain (b.2) or (b.3). Assume now that (2)
holds, if $C$ is an odd primitive cycle, by Lemma \ref{doublewheel}
we have (c.1), and if $C$ is an odd partial band, by Lemma
\ref{2oddring} we have (c.2).

$(\Leftarrow)$ It follows from Lemmas \ref{12sumas},
\ref{partialwheel},  \ref{doublewheel}, \ref{civertexband} and
\ref{2oddring}. \QED
\end{demo}

\bigskip Simis, Vasconcelos and Villarreal \cite{SVV} characterized the normality of $k[G]$ in the
following way.

\begin{Theorem}\label{ho} If $G$ is connected, then $k[G]$ is normal if
and only if every two vertex disjoint odd cycles are connected by an
edge.
\end{Theorem}

\medskip

From this description, one deduces that if a graph $G = [C; R]$
verifies that $k[G]$ is normal, then $C$ is connected. Moreover, $G$
cannot have a CI-vertex-band as induced subgraph.

\medskip

Thus, we can conclude the following results, which are the normal
versions of Theorem \ref{teoremaestructura} and Theorem
\ref{ultimo}.

\begin{Corollary}\label{teoremaestructuranormal}Let $G = [C; R]$ be a connected graph
such that $k[G]$ is normal. If $G$ is a complete intersection, then
\begin{itemize}
\item $R$ is a bipartite ring graph, and
\item $C$ is either the empty graph,
an odd primitive cycle or an odd partial band.
\end{itemize}
\end{Corollary}

\medskip

\begin{Corollary}\label{ultimonormal}Let $G = [C; R]$ be a connected graph such that $k[G]$ normal.
If $R$ is $2$-connected, then $G$ is a complete intersection if and
only if $G$ is one of the following graphs:
\begin{enumerate}
\item[{\rm (a)}] a bipartite ring graph,
\item[{\rm (b)}] a $1$-clique-sum of a bipartite ring graph and either
\begin{itemize} \item[{\rm (b.1)}] a CI-odd-partial-wheel, or \item[{\rm (b.2)}] a $1$-clique-sum of an odd partial band and an
edge,
\end{itemize}
\item[{\rm (c)}] a $2$-clique-sum of a bipartite ring graph and either \begin{itemize}  \item[{\rm (c.1)}] a CI-double-wheel, or
\item[{\rm (c.2)}] an odd partial band. \end{itemize}
\end{enumerate}
\end{Corollary}




\end{document}